\documentclass{amsproc}
\newtheorem{theorem}{Theorem}[section]
\newtheorem{corollary}[theorem]{Corollary}

\theoremstyle{definition}

\theoremstyle{remark}

\numberwithin{equation}{section}

\newcommand{\nwc}{\newcommand}
\nwc{\nwt}{\newtheorem}
\nwt{coro}{Corollary}
\nwt{ex}{Example}
\nwt{prop}{Proposition}
\nwt{defin}{Definition}

\nwc{\mf}{\mathbf} %Latex (as in \bf not tilted math letters)
\nwc{\blds}{\boldsymbol} %Latex 
\nwc{\ml}{\mathcal} %Latex

\nwc{\lam}{\lambda}
\nwc{\del}{\delta}
\nwc{\Del}{\Delta}
\nwc{\Lam}{\Lambda}
\nwc{\elll}{\ell}
\nwc{\IA}{\mathbb{A}} %algebraic
\nwc{\IB}{\mathbb{B}} %ball
\nwc{\IC}{\mathbb{C}} %complex
\nwc{\ID}{\mathbb{D}} %Dedekind
\nwc{\IE}{\mathbb{E}} %Euklides
\nwc{\IF}{\mathbb{F}} %finite field
\nwc{\IG}{\mathbb{G}} %Gauss
\nwc{\IH}{\mathbb{H}} %Hilbert\N-subgroup
\nwc{\IN}{\mathbb{N}} %natural
\nwc{\IP}{\mathbb{P}} %prime
\nwc{\IQ}{\mathbb{Q}} %rational
\nwc{\IR}{\mathbb{R}} %real
\nwc{\IS}{\mathbb{S}} %sphere
\nwc{\IT}{\mathbb{T}} %torus
\nwc{\IZ}{\mathbb{Z}} %integers

\def\bbleft{{\mathchoice {[\mskip-3mu {[}} {[\mskip-3mu {[}}{[\mskip-4mu {[}}{[\mskip-5mu {[}}}}
\def\bbright{{\mathchoice {]\mskip-3mu {]}} {]\mskip-3mu {]}}{]\mskip-4mu {]}}{]\mskip-5mu {]}}}}
\nwc{\setK}{\bbleft 1,K \bbright}
\nwc{\setN}{\bbleft 1,\cN \bbright}
 
\nwc{\va}{{\bf a}}
\nwc{\vb}{{\bf b}}
\nwc{\vc}{{\bf c}}
\nwc{\vd}{{\bf d}}
\nwc{\ve}{{\bf e}}
\nwc{\vf}{{\bf f}}
\nwc{\vg}{{\bf g}}
\nwc{\vh}{{\bf h}}
\nwc{\vi}{{\bf i}}
\nwc{\vI}{{\bf I}}
\nwc{\vj}{{\bf j}}
\nwc{\vk}{{\bf k}}
\nwc{\vl}{{\bf l}}
\nwc{\vm}{{\bf m}}
\nwc{\vM}{{\bf M}}
\nwc{\vn}{{\bf n}}
\nwc{\vo}{{\it o}}
\nwc{\vp}{{\bf p}}
\nwc{\vq}{{\bf q}}
\nwc{\vr}{{\bf r}}
\nwc{\vs}{{\bf s}}
\nwc{\vt}{{\bf t}}
\nwc{\vu}{{\bf u}}
\nwc{\vv}{{\bf v}}
\nwc{\vw}{{\bf w}}
\nwc{\vx}{{\bf x}}
\nwc{\vy}{{\bf y}}
\nwc{\vz}{{\bf z}}
\nwc{\bal}{\blds{\alpha}}
\nwc{\bep}{\blds{\epsilon}}
\nwc{\barbep}{\overline{\blds{\epsilon}}}
\nwc{\bnu}{\blds{\nu}}
\nwc{\bmu}{\blds{\mu}}
\nwc{\bet}{\blds{\eta}}

\nwc{\bk}{\blds{k}}
\nwc{\bm}{\blds{m}}
\nwc{\bM}{\blds{M}}
\nwc{\bp}{\blds{p}}
\nwc{\bq}{\blds{q}}
\nwc{\bn}{\blds{n}}
\nwc{\bv}{\blds{v}}
\nwc{\bw}{\blds{w}}
\nwc{\bx}{\blds{x}}
\nwc{\bxi}{\blds{\xi}}
\nwc{\by}{\blds{y}}
\nwc{\bz}{\blds{z}}

\nwc{\cA}{\ml{A}}
\nwc{\cB}{\ml{B}}
\nwc{\cC}{\ml{C}}
\nwc{\cD}{\ml{D}}
\nwc{\cE}{\ml{E}}
\nwc{\cF}{\ml{F}}
\nwc{\cG}{\ml{G}}
\nwc{\cH}{\ml{H}}
\nwc{\cI}{\ml{I}}
\nwc{\cJ}{\ml{J}}
\nwc{\cK}{\ml{K}}
\nwc{\cL}{\ml{L}}
\nwc{\cM}{\ml{M}}
\nwc{\cN}{\ml{N}}
\nwc{\cO}{\ml{O}}
\nwc{\cP}{\ml{P}}
\nwc{\cQ}{\ml{Q}}
\nwc{\cR}{\ml{R}}
\nwc{\cS}{\ml{S}}
\nwc{\cT}{\ml{T}}
\nwc{\cU}{\ml{U}}
\nwc{\cV}{\ml{V}}
\nwc{\cW}{\ml{W}}
\nwc{\cX}{\ml{X}}
\nwc{\cY}{\ml{Y}}
\nwc{\cZ}{\ml{Z}}

\nwc{\fA}{\mathfrak{a}}
\nwc{\fB}{\mathfrak{b}}
\nwc{\fC}{\mathfrak{c}}
\nwc{\fD}{\mathfrak{d}}
\nwc{\fE}{\mathfrak{e}}
\nwc{\fF}{\mathfrak{f}}
\nwc{\fG}{\mathfrak{g}}
\nwc{\fH}{\mathfrak{h}}
\nwc{\fI}{\mathfrak{i}}
\nwc{\fJ}{\mathfrak{j}}
\nwc{\fK}{\mathfrak{k}}
\nwc{\fL}{\mathfrak{l}}
\nwc{\fM}{\mathfrak{m}}
\nwc{\fN}{\mathfrak{n}}
\nwc{\fO}{\mathfrak{o}}
\nwc{\fP}{\mathfrak{p}}
\nwc{\fQ}{\mathfrak{q}}
\nwc{\fR}{\mathfrak{r}}
\nwc{\fS}{\mathfrak{s}}
\nwc{\fT}{\mathfrak{t}}
\nwc{\fU}{\mathfrak{u}}
\nwc{\fV}{\mathfrak{v}}
\nwc{\fW}{\mathfrak{w}}
\nwc{\fX}{\mathfrak{x}}
\nwc{\fY}{\mathfrak{y}}
\nwc{\fZ}{\mathfrak{z}}

\nwc{\tA}{\widetilde{A}}
\nwc{\tB}{\widetilde{B}}
\nwc{\tE}{E^{\vareps}}
\nwc{\tk}{\tilde k}
\nwc{\tN}{\tilde N}
\nwc{\tP}{\widetilde{P}}
\nwc{\tQ}{\widetilde{Q}}
\nwc{\tR}{\widetilde{R}}
\nwc{\tV}{\widetilde{V}}
\nwc{\tW}{\widetilde{W}}
\nwc{\ty}{\tilde y}
\nwc{\teta}{\tilde \eta}
\nwc{\tdelta}{\tilde \delta}
\nwc{\tlambda}{\tilde \lambda}
\nwc{\ttheta}{\tilde \theta}
\nwc{\tvartheta}{\tilde \vartheta}
\nwc{\tPhi}{\widetilde \Phi}
\nwc{\tpsi}{\tilde \psi}
\nwc{\tmu}{\tilde \mu}

\nwc{\To}{\longrightarrow} %limits

\nwc{\ad}{\rm ad}
\nwc{\eps}{\epsilon}
\nwc{\ep}{\epsilon}
\nwc{\vareps}{\varepsilon}

\def\ep{\epsilon}

\def\sq2{\sqrt{2}}

\def\t2{{\mathbb T}^2}
\def\s2{{\mathbb S}^2}

\def\T{\mathbb{T}}
\def\R{\mathbb{R}}

\def\Z{\mathbb{Z}}

\def\O{\mathcal{O}}

\nwc{\lap}{\bigtriangleup}
\nwc{\rest}{\restriction}
\nwc{\Diff}{\operatorname{Diff}}
\nwc{\diam}{\operatorname{diam}}
\nwc{\Res}{\operatorname{Res}}
\nwc{\Spec}{\operatorname{Spec}}
\nwc{\Vol}{\operatorname{Vol}}
\nwc{\Op}{\operatorname{Op}}
\nwc{\supp}{\operatorname{supp}}
\nwc{\Span}{\operatorname{span}}

\nwc{\dia}{\varepsilon}
\nwc{\cut}{f}
\nwc{\qm}{u_\hbar}

\def\hto0{\xrightarrow{\hbar\to 0}}

\def\rto0{\xrightarrow{r\to 0}}

\nwc{\la}{\langle}
\nwc{\ra}{\rangle}
\nwc{\lp}{\left(}
\nwc{\rp}{\right)}

\nwc{\bequ}{\begin{equation}}
\nwc{\be}{\begin{equation}}
\nwc{\ben}{\begin{equation*}}
\nwc{\bea}{\begin{eqnarray}}
\nwc{\bean}{\begin{eqnarray*}}
\nwc{\bit}{\begin{itemize}}
\nwc{\bver}{\begin{verbatim}}

\nwc{\eequ}{\end{equation}}
\nwc{\ee}{\end{equation}}
\nwc{\een}{\end{equation*}}
\nwc{\eea}{\end{eqnarray}}
\nwc{\eean}{\end{eqnarray*}}
\nwc{\eit}{\end{itemize}}
\nwc{\ever}{\end{verbatim}}

\begin{document}

\title[Semiclassical measures and the Schr\"odinger flow]{The dynamics of the Schr\"{o}dinger flow from the point of view of semiclassical measures}

\author{Nalini Anantharaman}
\address{Universit\'{e} Paris-Sud 11, Math\'{e}matiques, B\^{a}t. 425, 91405 ORSAY
CEDEX, FRANCE}
\email{Nalini.Anantharaman@math.u-psud.fr}
\thanks{N. Anantharaman wishes to
acknowledge the support of Agence Nationale de la Recherche, under
the grant ANR-09-JCJC-0099-01.}
\author{Fabricio Maci\`a}
\address{Universidad Polit\'{e}cnica de Madrid. DCAIN, ETSI Navales. Avda. Arco de la
Victoria s/n. 28040 MADRID, SPAIN} \email{Fabricio.Macia@upm.es}
\thanks{F. Maci{\`a} was supported by grants MTM2007-61755, MTM2010-16467 (MEC)}

%    General info
\date{}


\begin{abstract}
On a compact Riemannian manifold, we study the various dynamical
properties of the Schr\"odinger flow $(e^{it\Delta/2})$, through
the notion of {\em semiclassical measures} and the
quantum-classical correspondence between the Schr\"odinger
equation and the geodesic flow. More precisely, we are interested
in its high-frequency behavior, as well as its regularizing and
unique continuation-type properties. We survey a variety of
results illustrating the difference between positive, negative and
vanishing curvature.\end{abstract}

\maketitle

\section{Introduction\label{s:intro}}
Let $\left(  M,g\right)  $ be a smooth, $d$-dimensional, complete
manifold. Denote by $\Delta=\operatorname*{div}\left(
\nabla_{g}\cdot\right)  $ the Laplace-Beltrami operator and
consider the following linear Schr\"{o}dinger
equation on $M$:%
\begin{equation}
\left\{
\begin{array}
[c]{l}%
i\partial_{t}u\left(  t,x\right)  +\dfrac{1}{2}\Delta u\left(
t,x\right)
=0,\quad\left(  t,x\right)  \in\mathbb{R}\times M,\smallskip\\
u|_{t=0}=u^{0}\in L^{2}\left(  M\right)  .
\end{array}
\right.  \label{e:schrod}%
\end{equation}
Since $\left(  M,g\right)  $ is complete, $\Delta$ is an
essentially self-adjoint operator on $L^{2}\left(  M\right)  $ and
the initial value problem (\ref{e:schrod}) has a unique solution
$u\in C\left(  \mathbb{R} ;L^{2}\left(  M\right)  \right)  $. The
corresponding flow, the \emph{Schr\"{o}dinger flow}, is denoted by
$e^{it\Delta/2}$; recall that each operator $e^{it\Delta/2}$ is
unitary on $L^{2}\left(  M\right)  $, and in particular, for every
$t\in\mathbb{R}$ and $u^{0}\in L^{2}\left(  M\right) $,
\begin{equation}
\left\Vert e^{it\Delta/2}u^{0}\right\Vert _{L^{2}\left(  M\right)
}=\left\Vert u^{0}\right\Vert _{L^{2}\left(  M\right)  }.
\label{e:conservation}%
\end{equation}
When $M$ is compact a little more can be said: the solutions to
(\ref{e:schrod}) can be expressed in terms of eigenvalues and
eigenfunctions of $\Delta$ and the dynamics of $e^{it\Delta/2}$
turns out to be almost periodic~: if $(\varphi_j)_{j\in{\mathbb
N}}$ is an orthonormal basis formed of eigenfunctions on $\Delta$
(with $\Delta \varphi_j=-\lambda_j \varphi_j$), we can write
\begin{equation}
\label{e:develop} e^{it\Delta/2}u^{0} =\sum_{j\in {\mathbb
N}}e^{-it\lambda_j/2} \left( u^0 | \varphi_j \right) \varphi_j,
\end{equation}
where $\left(  \cdot|\cdot\right)  $ denotes the scalar product in
$L^{2}\left(  M\right)  $. However, this expression brings little
geometric information about the propagation properties of
$e^{it\Delta/2}$~: usually, the eigenfunctions  $\varphi_j$ are
not explicit, and even when they are (for instance, in the case of
a flat torus), oscillatory sums such as \eqref{e:develop} are
complicated objects.\footnote{In fact, identity \eqref{e:develop}
is often used to obtain information about the eigenfunctions from
the geometric description of the propagator $e^{it\Delta/2}$
\cite{NAQUE, AnNon, MaciaZoll}.} We shall return to the problem of
characterising the structure of eigenfunctions at the end of this
introduction and throughout the rest of this article; however, let
us state, from the very beginning, that this is not the point of
view that we are adopting (although, obviously, the dynamics of
the propagator $e^{it\Delta/2}$ and the properties of
eigenfunctions are closely related).

The issues we shall address here are aimed to obtain a better
understanding of the dynamics of $e^{it\Delta/2}$ and its relation
to the geometry of $\left( M,g\right)  $. In particular, we shall
be dealing with those aspects related to the high-frequency
behavior of $e^{it\Delta/2}$.

Let us describe precisely the main object of our study, before we
discuss in more detail the motivations that have guided us.
Consider
a sequence $\left( u_{n}^{0}\right)  $ of initial data in
$L^{2}\left(  M\right)  $ with $\left\Vert u_{n}^{0}\right\Vert
_{L^{2}\left(  M\right)  }=1$; we shall focus on the asymptotic
behavior, as $n\rightarrow\infty$, of the densities
\begin{equation}
\left\vert e^{it\Delta/2}u_{n}^{0}\right\vert ^{2}. \label{e:densities}%
\end{equation}
As we shall see, it is difficult to understand the behavior of
this quantity for individual $t$, but much more can be said if we
average w.r.t. $t$. Note that $\left\vert
e^{it\Delta/2}u_{n}^{0}\right\vert ^{2}\in L^{1}\left( M\right)  $
for every $t\in\mathbb{R}$, and, because of (\ref{e:conservation})
it can be identified to an element of $\mathcal{P}\left(  M\right)
$, the set of probability measures in $M$. Therefore, $\left(
\left\vert e^{it\Delta /2}u_{n}^{0}\right\vert ^{2}\right)  $ is a
sequence in $C\left( \mathbb{R};\mathcal{P}\left(  M\right)
\right)  $, and the Banach-Alaoglu
theorem ensures that it is compact in $L^{\infty}\left(  \mathbb{R}%
;\mathcal{M}\left(  M\right)  \right)  $ for the weak-$\ast$
topology.\footnote{Given a metric space $X$, we shall respectively
denote by $\mathcal{M}\left(  X\right)  $, $\mathcal{M}_{+}\left(
X\right)  $ and $\mathcal{P}\left(  X\right)  $ the set of Radon
measures, positive Radon measures and probability Radon measures
on $X$.}

In particular, there always exist a subsequence
$(u_{n^{\prime}}^{0})$ and a measure $\nu\in\linebreak
L^{\infty}\left(  \mathbb{R};\mathcal{M}_{+}\left( M\right)
\right)  $ such that
\begin{equation}
\int_{a}^{b}\int_{M}\chi\left(  x\right)  \left\vert e^{it\Delta
/2}u_{n^{\prime}}^{0}\right\vert ^{2}dt\rightarrow\int_{a}^{b}\int_{M}%
\chi\left(  x\right)  \nu\left(  t,dx\right)  dt,\quad\text{as
}n^{\prime
}\rightarrow\infty\text{,} \label{e:wL}%
\end{equation}
for every $\chi\in C_{c}\left(  M\right)  $ and
$a,b\in\mathbb{R}$. If $M$ is compact then $\nu\left(
t,\cdot\right)  $ is in fact a probability measure for a.e.
$t\in\mathbb{R}$. In general, the sequence $(u_{n^{\prime}}^{0})$
does not converge strongly in $L^{2}\left(  M\right)  $ and in
consequence, the measure $\nu$ may be singular with respect to the
Riemannian measure. The singular part of $\nu$ describes the
regions in $M$ on which the sequences of densities
(\ref{e:densities}) concentrates.\medskip\

Here we shall describe some results related to the question of
understanding the structure of the measures $\nu$ that arise in
this way. More precisely, we shall focus in aspects such as:

\begin{itemize}
\item The dependence of $\nu$ on the initial data $\left(
u_{n}^{0}\right) $. Is there a propagation law relating $\nu$ to
some limiting object obtained from the sequence $\left(
u_{n}^{0}\right)  $?

\item The regularity of $\nu$. Under which conditions on the
geometry of $\left(  M,g\right)  $ or on the structure of $\left(
u_{n}^{0}\right)  $ is it possible to ensure that the measure
$\nu$ is more regular than \emph{a priori }expected? For instance,
$\nu\in L^{p}\left(  \left[  a,b\right] \times M\right)  $ for
some $p>2$ and $a,b\in\mathbb{R}$.

\item The structure of the support of $\nu$. Which closed sets
$U\subset M$ can be the support of a measure $\nu$ obtained
through (\ref{e:wL}) for some sequence $\left(  u_{n}^{0}\right)
$?
\end{itemize}

One expects that the answer to these questions will strongly
depend on the geometry of $\left(  M,g\right)  $ and, in
particular, on the dynamics of the geodesic flow on the cotangent
bundle $T^{\ast}M$. Here we shall review some results obtained by
the authors in two different, and somewhat extremal, situations:
the cases of completely integrable (Sections \ref{s:zoll},
\ref{s:torus}) and Anosov geodesic flows (Section
\ref{s:negcurv}). From the point of view of manifolds of constant
sectional curvature, this corresponds to the cases of nonnegative
and negative sectional curvature, respectively. These results are
expressed in terms of semiclassical (or Wigner) measures, whose
main properties are recalled in Section \ref{s:ltsemiclassical}.

As we already mentioned, our motivation for addressing these
issues comes from the study of the dynamics of the linear
Schr\"{o}dinger equation, and more precisely on the following
three aspects.

\begin{enumerate}
\item The high-frequency dynamics of $e^{it\Delta/2}$ and its
relation to the quantum-classical correspondence principle and the
semiclassical limit of quantum mechanics.

\item The analysis of the dispersive properties of
$e^{it\Delta/2}$, and in particular the validity of Strichartz
estimates on a general Riemannian manifold.

\item The validity of observability or quantitative unique
continuation estimates for $e^{it\Delta/2}$.
\end{enumerate}

Before proceeding to describe the results, let us give a more
detailed description of each of these questions and, to conclude
this introduction, clarify how the problem addressed here is
related to other questions in Spectral Geometry that have been
widely studied in the literature (random initial data,
eigenfunction limits, and pair-correlation eigenvalue statistics).

\subsection{The quantum-classical correspondence principle and the
semiclassical limit}

The Schr\"{o}dinger equation (\ref{e:schrod}) is a mathematical
model for the propagation of a free quantum particle whose motion
is constrained to $M$. If $u$ is a solution to (\ref{e:schrod})
then for every measurable set $U\subset
M$ and every $t\in\mathbb{R}$, the quantity%
\begin{equation}\label{e:defproba}
\int_{U}\left\vert u\left(  t,x\right)  \right\vert ^{2}dx
\end{equation}
is the probability for the particle that was at $t=0$ at the state
$u^{0}$, to be in the region $U$ at time $t$. The
\emph{quantum-classical correspondence principle} asserts that if
the characteristic length of the oscillations of $u$ is very
small, then the dynamics of $\left\vert u\left(  t,\cdot\right)
\right\vert ^{2}$ can be deduced from that of the corresponding
classical system, that is, the geodesic flow $g^{t}$ on the
cotangent bundle $T^{\ast}M$ of $\left(  M,g\right)  $.

In order to develop a rigorous mathematical theory, we must
precise what we mean by a characteristic length of oscillations.
To do so, it is convenient to replace the initial datum $u^{0}$ by
a sequence $\left(  u_{n}^{0}\right)  $ of initial data with
$\left\Vert u_{n}^{0}\right\Vert _{L^{2}\left(  M\right) }=1$. Let
$\mathbf{1}_{\left[  0,1\right]  }$ denote the characteristic
function of the interval $\left[  0,1\right]  $. Chose a sequence
$\left( h_{n}\right)  $ of positive reals such that
$h_{n}\rightarrow0$ as
$n\rightarrow\infty$ and:%
\begin{equation}
\lim_{n\rightarrow\infty}\left\Vert \mathbf{1}_{\left[  0,1\right]  }%
(h_{n}\sqrt{-\Delta})u_{n}^{0}\right\Vert _{L^{2}\left(  M\right)
}=1;
\label{e:h-osc}%
\end{equation}
note that such a sequence always exists, by the spectral theorem
for self-adjoint operators on Hilbert space. If (\ref{e:h-osc})
holds we say that $\left(  h_{n}\right)  $ is a characteristic
length-scale for the oscillations of $\left(  u_{n}^{0}\right)  $,
or, following the terminology in \cite{GerardMesuresSemi91,
GerLeich93}, that $\left(  u_{n}^{0}\right)  $ is $\left(
h_{n}\right)  $\emph{-oscillating}.

As an example, let $\left(  x_{0},\xi_{0}\right)  \in T^{\ast}M$
and $u_{n}^{0}\in L^{2}\left(  M\right)  $ be supported on a
coordinate patch
around $x_{0}$ such that in coordinates:%
\begin{equation}
u_{n}^{0}\left(  x\right)  =\frac{1}{h_{n}^{d/4}}\rho\left(  \frac{x-x_{0}%
}{\sqrt{h_{n}}}\right)  e^{i\frac{\xi_{0}\cdot x}{h_{n}}},\label{e:WP}%
\end{equation}
where $\left(  h_{n}\right)  $ is a sequence of positive reals
tending to zero and $\rho$ is taken to have $\left\Vert
u_{n}^{0}\right\Vert _{L^{2}\left( M\right)  }=1$. The function
$u_{n}^{0}$ is usually called a wave-packet or coherent state
centered at $\left(  x_{0},\xi_{0}\right)  $. If $\Vert\xi
_{0}\Vert_{x_{0}}=1$,\footnote{In what follows, the expression
$\left\Vert \xi\right\Vert _{x}$ will denote the norm of $\left(
x,\xi\right)  \in T^{\ast}M$ induced by the Riemannian metric of
$\left(  M,g\right)  $.} then $\left(  u_{n}^{0}\right)  $ is
$\left(  h_{n}\right)  $-oscillating in the sense introduced
above. A manifestation of the correspondence principle is the
following classical result~: for any fixed $t\in\mathbb{R}$:%
\begin{equation}
\left\vert e^{ih_{n}t\Delta/2}u_{n}^{0}\right\vert
^{2}\rightharpoonup
\delta_{x\left(  t\right)  },\quad\text{as }n\rightarrow\infty\text{,}%
\label{e:semiclassicWP}%
\end{equation}
where $x\left(  t\right)  $ is the projection on $M$ of the orbit
$g^{t}\left(  x_{0},\xi_{0}\right)  $ of the geodesic flow.
Therefore, in the limit $n\rightarrow\infty$ the probability
densities $\left\vert e^{ih_{n}t\Delta/2}u_{n}^{0}\right\vert
^{2}$ become concentrated on the classical trajectory $x\left(
t\right)  $. Note that the time scale considered in the limit
(\ref{e:semiclassicWP}) is $h_{n}t$, which is proportional to the
characteristics length of oscillations of $\left( u_{n}^{0}\right)
$ and therefore tends to zero. An analogous result holds for more
general, $h_{n}$-oscillating sequences of initial data, see
Section \ref{s:ltsemiclassical} below; we shall refer to this as
the \emph{semiclassical limit}.

The convergence in (\ref{e:semiclassicWP}) is locally uniform in
$t\in\mathbb{R}$. Due to the dispersive nature of $e^{it\Delta/2}$
one cannot expect that (\ref{e:semiclassicWP}) holds uniformly in
time: for fixed $n$ and as $t$ increases, the wave-packet
$e^{ih_{n}t\Delta/2}u_{n}^{0}$ will become less and less
concentrated around $x\left(  t\right)  $.
The study of the simultaneous limits $h_{n}\To 0$ and $t\To \infty$ is a notoriously difficult problem. In the most general framework, it is known
\cite{CombescureRobertWP, BambusiGraffiPaul, HagedornJoye99,
HagedornJoye00, BouzouinaRobert} that (\ref{e:semiclassicWP})
holds uniformly for
\begin{equation}
\left\vert t\right\vert \leq T_{\text{E}}^{h_{n}}:=\left(
1-\delta\right)
\lambda_{\text{max}}^{-1}\log\left(  1/h_{n}\right)  ,\label{e:ehrenfest}%
\end{equation}
where $\delta\in\left(  0,1\right)  $ and $\lambda_{\text{max}}$
stands for the maximal expansion rate of the geodesic flow on the
spheres $\{\left\Vert \xi\right\Vert _{x}=\left\Vert
\xi_{0}\right\Vert _{x_{0}}\}$. This upper bound
$T_{\text{E}}^{h_{n}}$, known as the \emph{Ehrenfest time}, has
been shown to be optimal for some one-dimensional systems, see
\cite{deBievreRobert, Lablee}.

For the Euclidean space $\mathbb{R}^{d}$ or the torus
$\mathbb{T}^{d}$ equipped with the flat metric (or, more
generally, when the geodesic flow of $\left(  M,g\right)  $ is
completely integrable), it is possible to show that the
convergence in (\ref{e:semiclassicWP}) is uniform up to times
$\left\vert t\right\vert \leq Ch^{-1/2+\delta}$ for any
$\delta>0$, see \cite{BouzouinaRobert}. We stress the fact that
having an explicit expression of $e^{it\Delta/2}$ does not
necessarily make the study of the probability measures
\eqref{e:defproba} easier. To illustrate this phenomenon, let us
consider a gaussian coherent state on $\mathbb{R}^{d}$,
\begin{equation*}
v_{n}^{0}\left(  x\right)  =\frac{1}{h_{n}^{d/4}} e^{-\frac{|x-x_0|^2}{2h_n}}  e^{i\frac{\xi_{0}\cdot x}{h_{n}}}. %
\end{equation*}
One finds by an explicit calculation%
\begin{equation}
e^{ith_{n}\Delta/2}v_{n}^{0}\left(  x\right)  =\frac{1}{h_{n}^{d/4}%
(1+it)^{d/2}}e^{-\frac{\left\vert x-x_{0}-t\xi_{0}\right\vert ^{2}}%
{2h_{n}(1+t^{2})}}e^{i\frac{\phi\left(  t,x,x_{0},\xi_{0}\right)  }{2h_{n}}%
},\label{e:evolcoh}%
\end{equation}
with%
\[
\phi\left(  t,x,x_{0},\xi_{0}\right)  :=t\frac{\left\vert
x-x_{0}-t\xi _{0}\right\vert ^{2}}{\left(  1+t^{2}\right)
}-t\left\vert \xi_{0}\right\vert ^{2}+2\xi_{0}\cdot x.
\]
On the torus $\mathbb{T}^{d}=\mathbb{R}^{d}/2\pi\mathbb{Z}^{d}$,
this means that the evolution of the periodic wave packet
\[
u_{n}^{0}\left(  x\right)  =\frac{1}{h_{n}^{d/4}}\sum_{k\in2\pi\mathbb{Z}^{d}%
}e^{-\frac{|x-x_0-k|^{2}}{2h_{n}}}e^{i\frac{\xi_{0}\cdot(x-k)}{h_{n}}}%
\]
is given by the explicit expression
\begin{equation}
e^{ith_{n}\Delta/2}u_{n}^{0}\left(  x\right)  =\frac{1}{h_{n}^{d/4}%
(1+it)^{d/2}}\sum_{k\in2\pi\mathbb{Z}^{d}}e^{-\frac{\left\vert x-x_{0}-k%
-t\xi_{0}\right\vert ^{2}}{2h_{n}(1+t^{2})}}e^{i\frac{\phi\left(
t,x-k,x_{0},\xi_{0}\right)  }{2h_{n}}}.\label{e:propatorus}%
\end{equation}
For $\left\vert t\right\vert \leq Ch_{n}^{-1/2+\delta}$ it is
clear that the
associated probability measure \eqref{e:defproba} concentrates on the trajectory $x_{0}%
+t\xi_{0}$ (or on its image on the torus), but for $\left\vert
t\right\vert \geq h_{n}^{-1/2}$ these probability measures become
complicated objects due to the interferences between the different
terms in the sum \eqref{e:propatorus}. On compact negatively curved manifolds, there are also
examples of initial coherent or lagrangian states whose time
evolution is explicit up to times $t\sim h_{n}^{-2}$ \cite{Paul11,
Schub-largetimes}, but for which the associated probabilities
\eqref{e:defproba} are extremely complicated oscillatory sums.

One of our motivations is to study the probability measures
\eqref{e:defproba} at times $t$ for which the convergence
(\ref{e:semiclassicWP}) fails. Although this a very difficult
question for fixed $t$, it becomes more tractable if one performs
a time average. The problem we study consists in averaging the
probability measure \eqref{e:defproba} over a fixed time interval,
which means, with the semiclassical normalisation of time, to
average $\left\vert e^{ih_{n}t\Delta/2}u_{n}^{0}\right\vert ^{2}$
over time intervals of size $\sim h_n^{-1}$ or larger. In
particular, this averaging procedure allows us to go much beyond
the times where the individual $\left\vert
e^{ih_{n}t\Delta/2}u_{n}^{0}\right\vert ^{2}$ have been previously
studied.

We shall discuss and compare various geometries~: Zoll manifolds (Section \ref{s:zoll}), flat tori (Section \ref{s:torus}) and negatively curved manifolds (Section \ref{s:negcurv}).
Even when the geodesic
flow is completely integrable, important differences may occur, as
the analysis for the sphere $\mathbb{S}^{d}$ (or more generally,
of manifolds with periodic geodesic flow) and the torus
$\mathbb{T}^{d}$ shows, see Sections \ref{s:zoll} and
\ref{s:torus}.

\subsection{Dispersive properties of the Schr\"{o}dinger flow on a Riemannian
manifold}

By the word ``dispersion'', we mean that any solution to
the Schr\"{o}dinger equation (\ref{e:schrod}) can be expressed as
a superposition of waves propagating at different speeds,
depending on the characteristic frequencies of the initial datum.
For instance, when $M=\mathbb{R}^{d}$ any solution to
(\ref{e:schrod})
can be written as:%
\begin{equation}
e^{it\Delta/2}u_{0}\left(  x\right)  =\int_{\mathbb{R}^{d}}\widehat{u_{0}%
}\left(  \xi\right)  e^{i\xi\cdot\left(  x-t\frac{\xi}{2}\right)
}\frac{d\xi}{\left(
2\pi\right)  ^{d}}, \label{e:SchrodRd}%
\end{equation}
where $\widehat{u_{0}}$ stands for the Fourier transform of
$u_{0}$. This formula shows indeed that $e^{it\Delta/2}u_{0}$ is
built as a superposition of plane waves $e^{i\xi\cdot\left(
x-t\frac{\xi}{2}\right)  }$, travelling at velocity $\xi/2$. The
dispersion property is also seen very clearly in the expression
\eqref{e:evolcoh}, where we see that a coherent state initially
microlocalised around $(x_0, \xi_0)$ is less and less localized as
time evolves, while its $L^\infty$-norm decreases accordingly.

The representation formula (\ref{e:SchrodRd}) leads to the estimate:%
\begin{equation}
\left\Vert e^{it\Delta}u_{0}\right\Vert _{L^{\infty}\left(  \mathbb{R}%
^{d}\right)  }\leq\frac{C}{\left\vert t\right\vert
^{d/2}}\left\Vert u_{0}\right\Vert _{L^{1}\left(
\mathbb{R}^{d}\right)  }, \label{e:decay}
\end{equation}
which quantifies the
decay in time of solutions to (\ref{e:SchrodRd}) due to
dispersion. That estimate is in turn used to derive, by
interpolation with the conservation property
(\ref{e:conservation}), the commonly known as Strichartz
estimate:%
\begin{equation}
\left\Vert e^{it\Delta/2}u_{0}\right\Vert _{L^{p}\left(  \mathbb{R}_{t}%
\times\mathbb{R}_{x}^{d}\right)  }\leq C\left\Vert
u_{0}\right\Vert
_{L^{2}\left(  \mathbb{R}^{d}\right)  }, \label{e:StrichartzRd}%
\end{equation}
where
\begin{equation}
p=2\left(  1+\frac{2}{d}\right)  . \label{e:ExponentRd}%
\end{equation}
Estimate (\ref{e:StrichartzRd}) expresses that the singularities
(quantified by a Lebesgue norm) developed by a solution to the
Schr\"{o}dinger equation are better than what one would initially
expect based on the fact that $u_{0}\in L^{2}\left(
\mathbb{R}^{d}\right)  $. These estimates play a key role in the
well-posedness theory of semi-linear Schr\"{o}dinger equations,
see for instance \cite{GinibreNLS96, BourgainNLSBook,
CazenaveNLSBook, TaoDispersiveBook, GerardICM06} for an
introduction to this wide area of active research.

It is natural to wonder under which circumstances an estimate such
as (\ref{e:StrichartzRd}) holds if $\mathbb{R}^{d}$ is replaced by
a more general Riemannian manifold $\left(  M,g\right)  $. Or more
generally, how the geometry of $M$ affects the dispersive
character of the Schr\"{o}dinger flow. A first difficulty arises
in generalizing (\ref{e:StrichartzRd}) to a compact manifold~: as mentioned above, if $M$ is compact the
dynamics of $e^{it\Delta/2}$ turns out to be almost-periodic;
therefore, there is no hope for a global-in-time estimate to hold
in that case (clearly, no decay in time estimate as
(\ref{e:decay}) holds). But even if the time integral is replaced
by a local one, an estimate such as in (\ref{e:StrichartzRd}) may
still fail for any choice of $p>2$, as the example of the sphere
$\mathbb{S}^{2}$ shows, see \cite{BGT02}.

The validity of a Strichartz estimate:%
\begin{equation}
\left\Vert e^{it\Delta/2}u_{0}\right\Vert _{L^{p}\left(  \left[
0,1\right] \times M\right)  }\leq C\left\Vert u_{0}\right\Vert
_{L^{2}\left(  M\right)
}, \label{e:Strichartz}%
\end{equation}
on a compact manifold $M$ is related (in a somewhat loose manner)
to the regularity properties of the limit measures we introduced in
(\ref{e:wL}). Suppose that the Strichartz estimate
(\ref{e:Strichartz}) holds for some $p>2$, and let $\nu$ be a measure
obtained as in (\ref{e:wL}), for some sequence
$\left(  u_{n}^{0}\right)  $ of initial data in $L^{2}\left(
M\right)  $. Then the Strichartz inequality (\ref{e:Strichartz}) automatically implies that
$\nu\in L^{p/2}\left(  \left[  a,b\right]  \times M\right)  $, and
in particular, that $\nu$ is absolutely continuous with respect to
the Riemannian volume measure.
In other words, if one is able to construct a sequence $\left(  u_{n}%
^{0}\right)  $ that admits a measure $\nu$ as its limit
(\ref{e:wL}), such that $\nu$ has a non trivial singular part,
then this immediately shows that no estimate such as
(\ref{e:Strichartz}) holds for any $p>2$.

This is in fact the case when $\left(  M,g\right)  $ has periodic
geodesic flow (such a $\left(  M,g\right)  $ is called a Zoll
manifold), which proves that Strichartz estimates are false in
that case \cite{MaciaDispersion}. Note however that
frequency-dependent estimates (that is, with the $L^{2}\left(
M\right)  $-norm in the right-hand side of (\ref{e:Strichartz})
replaced by a Sobolev norm $H^{s}\left(  M\right)  $) still hold
in that case for exponents $s$ smaller than the one given by the
Sobolev embedding (see the works of Burq, G\'{e}rard and Tzvetkov
\cite{BGT04, BGT05, BGT05b}).

The situation is a bit different in the case of the flat torus $\mathbb{T}%
^{d}$. For $d=1$, a simple and elegant argument due to Zygmund
\cite{Zygmund74} shows that (\ref{e:Strichartz}) holds for $p=4$.
However, this is no longer the case for $d=1$, $p=6$ and $d=2$,
$p=4$ which are the exponents corresponding to the Euclidean space
(\ref{e:ExponentRd})$.$ Estimate (\ref{e:Strichartz}) fails in
those cases as shown by Bourgain \cite{BourgainStrichartz93,
BourgainIrrationalTori07}. Our results in that case
\cite{MaciaTorus, AnantharamanMacia}, developed in Section
\ref{s:torus}, imply that the measures obtained through
(\ref{e:wL}) are absolutely continuous with respect to the
Lebesgue measure, for any $d\geq1$; a proof of this fact based on
results on the distribution of lattice points on paraboloids is
indicated in the final remark of the article by Bourgain
\cite{BourgainQL97}. Moreover, in \cite{AnantharamanMacia} it is
shown that this absolute continuity result holds even for a more
general class of Hamiltonians defined on the flat torus.

\subsection{Observability and unique continuation for the Schr\"{o}dinger
flow\label{s:observability}}

A third aspect of the dynamics of the Schr\"{o}dinger flow, also
related to the properties of the limits (\ref{e:wL}), is the
validity of the \emph{observability }property, a quantitative
version of the unique continuation property that is relevant, for
instance, in Control Theory \cite{LionsSurvey88}, or Inverse
Problems \cite{IsakovBook}.

Let $T>0$ and $U\subset M$ be an open set; we say that the
Schr\"{o}dinger flow on $\left(  M,g\right)  $ satisfies the
observability property for $T$ and $U$ if a constant $C=C\left(
T,U\right)  >0$ exists, such that the
inequality%
\begin{equation}
\left\Vert u_{0}\right\Vert _{L^{2}\left(  M\right)  }^{2}\leq C\int_{0}%
^{T}\int_{U}\left\vert e^{it\Delta/2}u_{0}\left(  x\right)
\right\vert
^{2}dxdt \label{e:Observability}%
\end{equation}
holds for every initial datum $u_{0}\in L^{2}\left(  M\right)  $.
Clearly, the unique continuation property $\left[
e^{it\Delta/2}u_{0}|_{\left( 0,T\right)  \times U}\equiv0
\Longrightarrow u_{0}=0\right]  $ is a consequence of
(\ref{e:Observability}). However, (\ref{e:Observability}) also
implies a stronger stability property for the Schr\"{o}dinger
flow: two solutions to (\ref{e:schrod}) that are close to each
other in $\left( 0,T\right)  \times U$ (with respect to the
$L^{2}\left(  \left(  0,T\right) \times U\right)  $-norm) must
necessarily be issued from initial data that are also close in
$L^{2}\left(  \mathbb{R}^{d}\right)  $.

The following condition on $U$, sometimes referred to as the
\emph{Geometric Control Condition} is sufficient for the
observability property to hold for every $T>0$ , as shown by
Lebeau \cite{LebeauSchrod92} (see also
\cite{RauchTaylor, DGLControlNLS}).%
\begin{equation}%
\begin{array}
[c]{c}%
\text{There exists }L_{U}>0\text{ such that }\\
\text{every geodesic of }\left(  M,g\right)  \text{ of length
larger than
}L_{U}\text{ intersects }\overline{U}\text{.}%
\end{array}
\label{e:GCC}%
\end{equation}
In the particular case in which the geodesic flow of $\left(
M,g\right)  $ is periodic, it has been shown in
\cite{MaciaDispersion} that (\ref{e:GCC}) turns out to be also
necessary for observability.

However, this is not the case in general. For instance, when
$\left( M,g\right)  $ is the torus $\mathbb{T}^{d}$ equipped with
the flat canonical metric, a
result of Jaffard \cite{JaffardPlaques} (see also \cite{BurqZworskiBlackBox04}%
) shows that the observability property holds for every $T>0$ and
every open set $U\subset\mathbb{T}^{d}$, even if the Geometric
Control Condition fails. In the same direction, as proved in
\cite{AnRiv}, observability holds under conditions weaker than
(\ref{e:GCC}) for the case of manifolds of constant negative
sectional curvature, see Theorem \ref{t:observability}.

The observability property and the analysis of the limits
(\ref{e:wL}) are related as follows. Suppose that
(\ref{e:Observability}) holds for some $T$ and $U$. Any measure
$\nu$ obtained as a limit (\ref{e:wL}) would then satisfy
$\nu\left(  \left(  0,T\right)  \times U\right)  \geq1/C$. In
particular, the open set $\left(  0,T\right)  \times U$ must
intersect the support of every measure $\nu$ obtained by
(\ref{e:wL}) for any sequence of initial data $\left(
u_{n}^{0}\right)  $. Therefore, it is relevant in this context to
have detailed information on localization properties of the
measures obtained through (\ref{e:wL}).

%  In this r\'egime, we are interested in the propagation law
%satisfied by the limit measure and how it is related to the
%initial data $\left(  u_{n}^{0}\right)  $. As we shall see, such a
%propagation law is by no means universal;
% In this last case, we will see moreover that the
%propagation law obeyed by the limiting measures is not completely
%determined by the geodesic flow.
\subsection{Relations to other problems studied in the literature}
\subsubsection{Deterministic \emph{vs.} random sequences of initial data}
In this article we have addressed the problem of understanding how
the properties of the limits (\ref{e:wL}) depend on the geometry
of the ambient manifold. From this point of view, what we seek is
to prove properties of the limits $\nu$ that hold for any {\em
arbitrary} sequence $(u_n^0)$, thus reflecting the geometric
features of the propagator $e^{it\Delta/2}$.

In the present context {\em arbitrary} is not synonymous to {\em
random}. By the term {\em random} one can mean random sequences of
initial data, as in \cite{ZelditchQESphere} where it is shown that
for almost all orthonormal base of eigenfunctions of the laplacian
on the sphere, the limit (\ref{e:wL}) coincides with the standard
Riemannian volume. \emph{Random} can also refer to the fact that
the coefficients $(u^0|\varphi_j)$ in \eqref{e:develop} are random
variables~: for instance, independent centered gaussians. If one
is interested in the high-frequency r\'egime, one should restrict
to $\lambda_j$ in some interval $[E-\delta E, E+\delta E]$ and
take the limit $E\To +\infty$. If $\delta E\gg E^{1/2}$, it is
easy to show that, for any given $t$ and almost surely, $|u(t,
x)|^2 dx$ converges to the uniform measure on $M$ as $E\To
+\infty$. Note that this type of result is independent of the
geometry of $M$, and its scope is different from the type of
problem previously described here, namely that of characterising
the limits (\ref{e:wL}) for \emph{every} possible bounded
sequence. As particular cases, we deal with coherent states, Dirac
states, and the eigenfunctions of the laplacian $\varphi_j$
themselves.

\subsubsection{The case of eigenfunctions}
When the initial condition $u^0$ is an eigenfunction of the
laplacian $\varphi_j$, the probability measure \eqref{e:densities}
does not depend on $t$ and simply reads $|\varphi_j(x)|^2 dx$. The
behavior of these measures as $j\To +\infty$ has been the center
of much attention recently, in particular in the context of the
Quantum Unique Ergodicity conjecture (see Section \ref{s:negcurv}
and the references therein, as well as Sections \ref{s:zoll} and
\ref{s:torus} for results and references on the completely
integrable case). The study of the limit \eqref{e:densities} is
more general, and thus all the results we mention below also apply
to eigenfunctions.

We stress again the fact that having an explicit expression of the
eigenfunctions does not necessarily make the problem easy~:
consider for instance the case of the flat torus
$\T^d=\R^d/2\pi\Z^d$. An eigenfunction that satisfies $\Delta
\varphi=-\lambda \varphi$ can be decomposed as
\begin{equation}\label{e:torusef}\varphi(x)=\sum_{k\in \Z^d, |k|^2=\lambda}c_k e^{ik\cdot x}.\end{equation}
The spectral degeneracy of $\lambda$, that is, the number of
integral solutions of $|k|^2=\lambda$, gets unbounded as $\lambda$
grows. Sums of the form \eqref{e:torusef} and the corresponding
squares $|\varphi(x)|^2$ have been studied in
\cite{JakobsonTori97, MarklofSquares}.

\subsubsection{Level spacings and pair-correlation statistics}
Let $\operatorname{sp}\left(  -\Delta\right)  $ denote the
spectrum of the Laplace-Beltrami operator and for
$\lambda\in\operatorname{sp}\left( -\Delta\right)  $, write
$P_{\lambda}$ to denote the orthogonal projection from
$L^{2}\left(  M\right)  $ onto the eigenspace associated to
$\lambda$. As before, denote by $\left(  \lambda_{j}\right)
_{j\in\mathbb{N}}$ the eigenvalues of $-\Delta$ counted with their
multiplicities and by $\left( \varphi_{j}\right)
_{j\in\mathbb{N}}$ an orthonormal basis consisting of
eigenfunctions indexed accordingly. Using \eqref{e:develop}, we
see that for any $\theta\in L^{1}(\R)$ and $\chi\in C(M)$, the
expression
\begin{equation}
\int_{\R}\theta(t)\int_{M}\chi(x)|e^{it\Delta/2}u|^{2}(x)dx\label{e:timedep}%
\end{equation}
can be expanded into~:%
\begin{equation}
\sum_{\lambda,\lambda^{\prime}\in\operatorname{sp}\left(
-\Delta\right)
}\hat{\theta}\left(  \frac{\lambda-\lambda^{\prime}}{2}\right)  \int_{M}%
\chi(x)P_{\lambda}u\left(  x\right)
\overline{P_{\lambda^{\prime}}u\left(
x\right)  }dx,\label{e:samebis}%
\end{equation}
or, equivalently,
\begin{equation}
\sum_{j,j^{\prime}\in{\mathbb{N}}}\hat{\theta}\left(
\frac{\lambda _{j}-\lambda_{j^{\prime}}}{2}\right)  \left(
u|\varphi_{j}\right) \overline{\left(
u|\varphi_{j^{\prime}}\right)  }\int_{M}\chi(x)\varphi
_{j}(x)\overline{\varphi_{j^{\prime}}(x)}dx,\label{e:same}%
\end{equation}
where $\hat{\theta}$ denotes the Fourier transform of $\theta$. In
particular, if $\hat{\theta}$ is compactly supported, this
restricts our sum to bounded $\lambda-\lambda^{\prime}$ (resp.
$\lambda_{j}-\lambda_{j^{\prime}}$), and two natural questions
arise~:

\begin{enumerate}
\item[(1)] Can the study of (\ref{e:timedep}), (\ref{e:samebis})
be reduced to that of the matrix elements
$\int_{M}\chi(x)P_{\lambda}u\left(  x\right)
\overline{P_{\lambda^{\prime}}u\left(  x\right)  }dx$ (for bounded
$\lambda-\lambda^{\prime}$)~?

\item[(2)] Does the knowledge of the distribution of the pair
correlations $\lambda_{j}-\lambda_{j^{\prime}}$ help to gain some
insight in (\ref{e:timedep}), (\ref{e:same})~?
\end{enumerate}
To answer (1), it is quite clear from \eqref{e:samebis} that the
study of \eqref{e:timedep} amounts in some sense to a study of
matrix elements; However, it is not necessarily easier to study
the matrix elements than to study the time-dependent equation. We can note that, in the (very special) cases where the minimal spacing $\inf\left\{
\lambda-\mu:\lambda,\mu\in\operatorname{sp}\left( -\Delta\right)
,\lambda\neq\mu\right\}  $ is strictly positive (as is the case of
the sphere or the flat torus, for instance), \eqref{e:samebis}
takes a
particularly simple form:%
\[
\int_{\mathbb{R}}\theta\left(  t\right)  dt\sum_{\lambda\in\operatorname{sp}%
\left(  -\Delta\right)  }\int_{M}\chi\left(  x\right)  \left\vert
P_{\lambda }u\left(  x\right)  \right\vert ^{2}dx.
\]
This property has been exploited in \cite{MaciaZoll}, on certain
classes of Zoll manifolds, to characterise the accumulation points
of sequences of the form $\left\vert P_{\lambda}u\left(  x\right)
\right\vert ^{2}dx$ from the knowledge of the structure of the
limits of (\ref{e:timedep}) (this can also be deduced from the
fine study of the structure of $P_{\lambda}$ performed in
\cite{ZelditchZoll}). Also in \cite{MaciaZoll} it is shown how the
study of (\ref{e:timedep}) can be used to obtain information on
the off-diagonal matrix elements (\ref{e:samebis}) in the case of
Zoll manifolds. The relations between time averaging and
eigenvalue level spacing are further explored in the forthcoming
article \cite{AFM} in the context of completely integrable
systems.

In answer to (2), we first recall that the pair correlation
distribution is conjectured to be Poissonian in the completely
integrable case \cite{BerTab}; this has been proved in a certain
number of cases \cite{SarnakSpacings, VdK1, VdK2, VdK3,
MarklofRandomWalk, MarklofPoisson, MarklofPoisson2, EMM}. At the
opposite end of \textquotedblleft chaotic
systems\textquotedblright, e.g. the case of the laplacian on
negatively curved surfaces, the pair correlation distribution is
conjectured to be given by Random Matrix Theory \cite{RMatrix,
HO84}; there is no mathematical proof of this fact, but this is a
field of active current research \cite{Bogo, Sieb02,
SieberRichter}. In any case, we do not think that the pair
correlation distribution bears any obvious relevance to the
understanding of \eqref{e:timedep} or \eqref{e:same}. As we said,
the study of \eqref{e:timedep} is already of interest when $u$ is
itself an eigenfunction $\varphi_{j}$, in which case
\eqref{e:same} is just $\int _{M}\chi(x)|\varphi_{j}(x)|^{2}dx$
and the pair correlations play absolutely no role. To compare both
problems, we may also comment on the case of flat tori~: the pair
correlation problem for rational and irrational tori is different,
whereas the \textquotedblleft shape\textquotedblright\ of the
torus does not seem to play any r\^{o}le in our study of
\eqref{e:timedep}. We can also add, in the case of negatively
curved manifolds, that our knowledge of the pair correlation
distribution is purely conjectural, whereas we do have results
about \eqref{e:timedep} (\S \ref{s:negcurv}). Thus, the link
between the two problems that seems to arise when rewriting
\eqref{e:timedep} in the form \eqref{e:same} is only apparent. We believe
that the study of the pair correlation problem is even more
difficult.

\section{Semiclassical measures\label{s:ltsemiclassical}}
Consider again the example of a wave-packet sequence of initial
data $\left( u_{n}^{0}\right)  $ centered at point $\left(
x_{0},\xi_{0}\right)  $ in the
cotangent bundle $T^{\ast}M$ as defined in (\ref{e:WP}). In this case,%
\[
\left\vert u_{n}^{0}\right\vert ^{2}\rightharpoonup\delta_{x_{0}}%
,\quad\text{as }n\rightarrow\infty\text{;}%
\]
note that this holds independently of the direction of oscillation
$\xi_{0}\in T_{x_{0}}^{\ast}M$. However, as we saw
(\ref{e:semiclassicWP}), the densities $\left\vert
e^{ih_{n}t\Delta/2}u_{n}^{0}\right\vert ^{2}$ corresponding to the
evolution concentrate on the point $x\left(  t\right)  $ of the
geodesic of $M$ issued from $\left(  x_{0},\xi_{0}\right)  $.
Therefore, their limit does effectively depend on $\xi_{0}$. This
shows that there is no propagation law
relating the limit of the densities $\left\vert e^{it\Delta/2}u_{n}%
^{0}\right\vert ^{2}$ to that of the initial densities $\left\vert u_{n}%
^{0}\right\vert ^{2}$.

This difficulty is overcome by lifting the measure $\left\vert u_{n}%
^{0}\right\vert ^{2}dx$ to phase space $T^{\ast}M$, which allows
to keep track of the characteristic oscillation frequencies of
$u_{n}^{0}$. There are different procedures to accomplish this,
but all of them are equivalent for our purposes. Here we shall
focus on the one based on the \emph{Weyl quantization }(see
\cite{FollandPhaseSpace} for a comprehensive introduction). Let us
first discuss the definition in the case $M=\mathbb{R}^{d}$ and
then give the general case.

Starting from a function (a classical observable) $a\in
C_{c}^{\infty}\left( T^{\ast}\mathbb{R}^{d}\right)  $, the Weyl
quantization associates to $a$ the operators
$\operatorname*{Op}_{h}\left(  a\right)  $ (with the
\textquotedblleft semiclassical\textquotedblright\ parameter
$h>0$), that act
on tempered distributions $u\in\mathcal{S}^{\prime}\left(  T^{\ast}%
\mathbb{R}^{d}\right)  $ as follows:%
\[
\operatorname*{Op}\nolimits_{h}\left(  a\right)  u\left(  x\right)
:=\int_{\mathbb{R}^{d}}\int_{\mathbb{R}^{d}}a\left(  \frac{x+y}{2}%
,h\xi\right)  u\left(  y\right)  e^{i\xi\cdot\left(  x-y\right)
}dy\frac {d\xi}{\left(  2\pi\right)  ^{d}}.
\]
Those operators are uniformly bounded from $L^{2}\left(  \mathbb{R}%
^{d}\right)  $ into itself, in fact (see \cite{GerLeich93}):
\begin{equation}
\left\Vert \operatorname*{Op}\nolimits_{h}\left(  a\right)
\right\Vert _{\mathcal{L}\left(  L^{2}\left(  M\right)  \right)
}\leq C_{d}\left\Vert
a\right\Vert _{C^{d+1}\left(  T^{*}M\right)  }.\label{e:ophbdd}%
\end{equation}
Note that when $a$ only depends on the variable $x$, the
corresponding Weyl operator acts on functions by multiplication by
$a$. On the other hand, if $a$ only depends on $\xi$, then
$\operatorname*{Op}\nolimits_{h}\left(  a\right)
$ is simply the Fourier multiplier:%
\[
a\left(  hD_{x}\right)  u\left(  x\right)
=\int_{\mathbb{R}^{d}}a\left( h\xi\right)  \widehat{u}\left(
\xi\right)  e^{i\xi\cdot x}\frac{d\xi}{\left( 2\pi\right)  ^{d}}.
\]
One can extend this definition to functions $a\in
C_{c}^{\infty}\left( T^{\ast}M\right)  $ for a general manifold
$M$ by means of local coordinates and partitions of unity, see for
instance \cite{GerLeich93, EvansZworski}.

In what follows, $\left(  u_{n}^{0}\right)  $ will be a bounded
sequence in $L^{2}\left(  M\right)  $ and $\left(  h_{n}\right)  $
a sequence of positive reals tending do zero such that the
$h_{n}$-oscillation condition (\ref{e:h-osc}) is fulfilled. We
shall define a distribution $w_{h_{n}}$ on
$T^{\ast}M$, which is a \emph{lift} of the measure $\left\vert u_{n}%
^{0}\right\vert ^{2}dx$, in the sense that it projects down to
$\left\vert u_{n}^{0}\right\vert ^{2}dx$ under the canonical
projection $T^{\ast }M\longrightarrow M$. The action of the
distribution $w_{h_{n}}\in \mathcal{D}^{\prime}\left(
T^{\ast}M\right)  $ on a test function $a\in
C_{c}^{\infty}\left(  T^{\ast}M\right)  $ is given by:%
\[
\left\langle w_{h_{n}},a\right\rangle :=\left(  \operatorname*{Op}%
\nolimits_{h_{n}}\left(  a\right)  u_{n}^{0}|u_{n}^{0}\right)  .
\]
Usually, because of E.P. Wigner's seminal work \cite{Wigner32},
$w_{h_{n}}$ is called the \emph{Wigner distribution }of the
function $u_{n}^{0}$. The sequence of distributions $\left(
w_{h_{n}}\right) $ is uniformly bounded, as a consequence of
(\ref{e:ophbdd}). It turns out that any accumulation point of
$\left(  w_{h_{n}}\right)  $ (in the weak
topology of distributions) is a positive measure $\mu_{0}\in\mathcal{M}%
_{+}\left(  T^{\ast}M\right)  $ despite the fact that the
$w_{h_{n}}$ are not positive. See \cite{CdV85, Zel87,
GerardMesuresSemi91, LionsPaul, EvansZworski} for different proofs
of this non trivial result. Moreover, if some subsequence of
$(\left\vert u_{n}^{0}\right\vert ^{2})$ and $\left(
w_{h_{n}}\right)  $ converges respectively to some measures
$\nu_{0}\in\mathcal{M}_{+}\left(
M\right)  $ and $\mu_{0}\in\mathcal{M}_{+}\left(  T^{\ast}M\right)  $ then:%
\begin{equation}
\nu_{0}\left(  x\right)  =\int_{T_{x}^{\ast}M}\mu_{0}\left(
x,d\xi\right)  .
\label{e:proj}%
\end{equation}
This means that $\mu_{0}$ is also a lift of $\nu_{0}$. Usually,
$\mu_{0}$ is
called a \emph{semiclassical measure }of the sequence $\left(  u_{n}%
^{0}\right)  $. \emph{A priori}, there can be several
semiclassical measures, as different subsequences may have
different limits. It is easy to show that a wave-packet
(\ref{e:WP}) has a unique semiclassical measure which is
$\delta_{\left(  x_{0},\xi_{0}\right)  }$.

We now turn to the problem of computing semiclassical measures of
sequences of trajectories to the Schr\"{o}dinger flow. We shall
denote by $w_{h_{n}}\left( t\right)  $ the Wigner distribution of
$e^{it\Delta/2}u_{n}^{0}$. The main tool in this context is
Egorov's theorem, which relates the Schr\"{o}dinger group
$e^{it\Delta/2}$ to the geodesic flow $g^{t}$ (see
\cite{EvansZworski} for a proof).

\begin{theorem}
For every $a\in C_{c}^{\infty}\left(  T^{\ast}M\right)  $ there
exists a family $R_{h}\left(  t\right)  $ of bounded operators on
$L^{2}\left(
M\right)  $ such that%
\begin{equation}
e^{-ith\Delta/2}\operatorname*{Op}\nolimits_{h}\left(  a\right)
e^{ith\Delta/2 }=\operatorname*{Op}\nolimits_{h}\left(  a\circ
g^{t}\right)
+R_{h}\left(  t\right)  , \label{e:egorov}%
\end{equation}
and $\left\Vert R_{h}\left(  t\right)  \right\Vert
_{\mathcal{L}\left( L^{2}\left(  M\right)  \right) }\leq\rho\left(
\left\vert t\right\vert
\right)  h$ for some non-negative continuous function $\rho:\mathbb{R}%
_{+}\mathbb{\rightarrow R}_{+}$.
\end{theorem}

With this result at our disposal, it is not hard to derive a
propagation law for the time-scaled Wigner distributions
$w_{h_{n}}\left(  h_{n}t\right)  $.

\begin{theorem}
\label{t:semiclassicallimit}Let $\left(  u_{n}^{0}\right)  $ and
$\left( h_{n}\right)  $ be as above. It is possible to extract a
subsequence such
that, for every $t\in\mathbb{R}$,%
\[
w_{h_{n^{\prime}}}\left(  h_{n^{\prime}}t\right)
\rightharpoonup\mu
_{t}^{\text{\emph{sc}}}\text{,\quad as }n^{\prime}\rightarrow\infty\text{,}%
\]
where $\mu_{t}^{\text{\emph{sc}}}$ is a continuous family of
positive measures
in $\mathcal{M}_{+}\left(  T^{\ast}M\right)  $. Moreover, $\mu_{t}%
^{\text{\emph{sc}}}$ is transported along the geodesic flow of
$\left(
M,g\right)  $:%
\begin{equation}
\mu_{t}^{\text{\emph{sc}}}=\left(  g^{t}\right)
_{\ast}\mu_{0},\quad\text{for
every }t\in\mathbb{R}\text{.} \label{e:transport}%
\end{equation}

\end{theorem}

Identity (\ref{e:transport}) means that $\mu_{t}^{\text{sc}}$ is
obtained as
the push-forward of $\mu_{0}$ by the geodesic flow:%
\[
\int_{T^{\ast}M}a\left(  x,\xi\right)  \mu_{t}^{\text{sc}}\left(
dx,d\xi\right)  =\int_{T^{\ast}M}a\circ g^{t}\left(  x,\xi\right)
\mu _{0}\left(  dx,d\xi\right)  ,
\]
for every $a\in C_{c}\left(  T^{\ast}M\right)  $.

As a consequence of this result, the semiclassical measure of the
evolution of
a wave-packet initial datum (\ref{e:WP}) is $\delta_{g^{t}\left(  x_{0}%
,\xi_{0}\right)  }$. Using the projection identity (\ref{e:proj})
we deduce the propagation law (\ref{e:semiclassicWP}) stated in
the introduction.

More generally, Theorem \ref{t:semiclassicallimit} can be used to
compute the limit of $\left\vert
e^{ith_{n}\Delta/2}u_{n}^{0}\right\vert ^{2}$ but does not apply
to obtain that of $\left\vert e^{it\Delta/2}u_{n}^{0}\right\vert
^{2}$. This is related to the sensitivity to time dependence of
Egorov's theorem. The identity that is relevant to our analysis is
obtained by rescaling time of a factor $1/h_{n}$ in
(\ref{e:egorov}). The remainder $R_{h}\left(  t\right)  $ is only
known to go to zero as $h\rightarrow0^{+}$ uniformly for
$\left\vert t\right\vert \leq T_{\text{E}}^{h_{n}}$, where
$T_{\text{E}}^{h_{n}}$ is the Ehrenfest time defined in
(\ref{e:ehrenfest}), see \cite{BouzouinaRobert}. Therefore, it is
not possible to ensure that $R_{h}\left(  t/h_{n}\right)  $ will
tend to zero as $h\rightarrow0^{+}$. But even if (\ref{e:egorov})
is exact (\emph{i.e. }$R_{h}\left(  t\right) \equiv0$, as is the
case when $M=\mathbb{T}^{d}$), it is not easy to deal with
the operators $\operatorname*{Op}\nolimits_{h_{n}}\left(  a\circ g^{t/h_{n}%
}\right)  $, due to the fact that the functions $a\circ
g^{t/h_{n}}$ depend on $h_{n}$ and vary very rapidly as $h_{n}$
goes to zero.

This problem has been widely studied when $M$ is compact and the
initial data
are normalized eigenfunctions of the Laplacian:%
\[
-\Delta u_{n}^{0}=\lambda_{n}u_{n}^{0},\text{\qquad}\left\Vert u_{n}%
^{0}\right\Vert _{L^{2}\left(  M\right)  }^{2}=1.
\]
corresponding to a sequence of eigenvalues $\left(
\lambda_{n}\right)  $ that tends to infinity as
$n\rightarrow\infty$. In this case, because of (\ref{e:h-osc}), it
is natural to set $h_{n}:=\lambda_{n}^{-1/2}$ and it turns
out that, for every $t\in\mathbb{R}$:%
\[
\left\vert e^{it\Delta/2}u_{n}^{0}\right\vert ^{2}=\left\vert u_{n}%
^{0}\right\vert ^{2}\text{,\qquad}w_{h_{n}}\left(  t\right)
=w_{h_{n}}\left( 0\right)  .
\]
Since these quantities do not depend on $t$, Theorem
\ref{t:semiclassicallimit} shows that any semiclassical measure
$\mu_{0}$ of $\left(  u_{n}^{0}\right)  $ is invariant by the
geodesic flow: $\left( g^{s}\right)  _{\ast}\mu_{0}=\mu_{0}$ for
every $s\in\mathbb{R}$. Moreover, it can easily be proved that
$\mu_{0}$ is supported on the cosphere bundle $S^{\ast}M:=\left\{
\left(  x,\xi\right)  \in T^{\ast}M\;:\;\left\Vert \xi\right\Vert
_{x}^{2}=1\right\}  $; therefore we can view $\mu_{0}$ as an
element of $\mathcal{P}\left(  S^{\ast}M,g^{t}\right)  $, the set of $g^{t}%
$-invariant probability measures on $S^{\ast}M$.

The problem of identifying those measures in $\mathcal{P}\left(
S^{\ast }M,g^{t}\right)  $ that arise as semiclassical measures of
some sequence of eigenfunctions of the Laplacian has proven to be
very hard in general. In contrast with the propagation law
(\ref{e:transport}), global aspects of the dynamics of $g^{t}$
play a role; in particular, the global geometry of $\left(
M,g\right)  $ is relevant for this problem. Some results on this
issue will be reviewed later on in this article.

Let us now turn to the case of arbitrary initial data in
\eqref{e:densities}. A first difficulty one encounters when
dealing with the Wigner distributions $w_{h_{n}}\left(  t\right) $
of $e^{it\Delta/2}u_{n}^{0}$ is that, due to the highly
oscillating nature of the propagator $e^{it\Delta/2}$, it is in
general not possible to extract a subsequence such that
$w_{h_{n}}\left(  t\right)  $ converges for every
$t\in\mathbb{R}$. This difficulty can be overcome by viewing
$\left(  w_{h_{n}}\right)  $ as a sequence in $L^{\infty}\left(
\mathbb{R};\mathcal{D}^{\prime}\left(  T^{\ast}M\right)  \right) $
and considering its accumulation points with respect to the
weak-$\ast$ topology in that space. This is nothing else but
considering time averages of $w_{h_{n}}$. The following result
holds (see \cite{MaciaAv}).

\begin{theorem}
\label{t:convergence}Let $\left(  u_{n}^{0}\right)  $ and $\left(
h_{n}\right)  $ be as above. Then there exist a subsequence
$\left( u_{n^{\prime}}^{0}\right)  $ and a measure $\mu\in
L^{\infty}\left( \mathbb{R};\mathcal{M}_{+}\left( T^{\ast}M\right)
\right)  $ such that, for
every $\varphi\in L^{1}\left(  \mathbb{R}\right)  $,%
\[
\int_{\mathbb{R}}\varphi\left(  t\right)  w_{h_{n^{\prime}}}\left(
t\right) dt\rightharpoonup\int_{\mathbb{R}}\varphi\left(  t\right)
\mu(t)dt,\quad
\text{as }n\rightarrow\infty\text{.}%
\]
Moreover, for a.e. $s\in\mathbb{R}$ the measure $\mu(s)$ is
invariant by the geodesic flow.
\end{theorem}

The measure $\mu\in L^{\infty}\left(
\mathbb{R};\mathcal{M}_{+}\left( T^{\ast}M\right)  \right)  $ will
be called a \emph{semiclassical measure} associated to
$(e^{it\Delta/2}u_{n}^{0})$ (or \emph{time-dependent semiclassical
measure} when we want to stress the difference with the previous
definition). In full generality, there is no propagation law
relating the semiclassical measures of the initial data
$(u_{n}^{0})$ and the semiclassical measures of
$(e^{it\Delta/2}u_{n}^{0})$, as we shall see in \S \ref{s:torus}.

It should be noted that Theorem \ref{t:convergence} still holds
for the time rescaled distributions $w_{h_{n}}\left(
h_{n}\alpha_{n}t\right)  $ where $\left(  \alpha_{n}\right)  $ is
any sequence that tends to infinity as $n\rightarrow\infty$ (see
\cite{MaciaAv}). Here we restrict ourselves to the case
$\alpha_{n}=1/h_{n}$. The rest of the article is devoted to
understanding how the structure of the measures $\mu(t)$ depends
on the geometry of $\left( M,g\right)  $.

\section{Zoll manifolds\label{s:zoll}}
In this section we shall deal with manifolds whose geodesic flow
has the simplest possible dynamics. We shall assume that $\left(
M,g\right)  $ is a compact manifold all of whose geodesics are
closed. These are called Zoll manifolds, and the book
\cite{BesseZoll} provides a comprehensive treatment of these
geometries. It is known that the geodesic flow of a Zoll manifold
is periodic. Every manifold of positive constant sectional
curvature (that is, the sphere $\mathbb{S}^{d}$ and its quotients
\cite{WolfBook}) is a Zoll manifold; the same holds for compact
symmetric spaces of rank one, as the complex projective spaces. O.
Zoll constructed a real analytic Riemannian metric on the sphere
$\mathbb{S}^{2}$, which is not isometric to the canonical one, but
still has the property that every geodesic is closed. It should be
noted that the geodesic flow on the cotangent bundle of a Zoll
manifold is a completely integrable Hamiltonian system
\cite{Duran}.

The first result on the structure of the set of semiclassical
measures for solutions to the Schr\"{o}dinger equation is due to
Jakobson and Zelditch \cite{JakobsonZelditch}. These authors
consider the case $M=\mathbb{S}^{d}$, equipped with its canonical
metric, and study semiclassical measures arising from sequences of
eigenfunctions of the Laplacian~: they show that any
invariant probability measure in $\mathcal{P}\left(  S^{\ast}\mathbb{S}%
^{d},g^{t}\right)  $ can be obtained as the semiclassical measure
of some sequence of eigenfunctions. At the origin of the proof of
this result, there is the easy remark that the restriction to the
sphere of the following
harmonic polynomials:%
\[
\psi_{n}\left(  x\right)  =C_{n}\left(  x_{1}+ix_{2}\right)  ^{n}%
,\qquad\left\Vert \psi_{n}\right\Vert _{L^{2}\left(
\mathbb{S}^{d}\right) }=1,
\]
($n\in{\mathbb{Z}}$, $|n|\longrightarrow\infty$) concentrates on
the maximal circle $x_{1}^{2}+x_{2}^{2}=1$ on $\mathbb{S}^{d}$.
The semiclassical measure of $\left(  \psi_{n}\right)  $ is
concentrated on one of the two orbits of the geodesic flow, that
lie above the aforementioned geodesic in the unit cotangent bundle
(the orientation depends on the sign of $n$). Thus, any invariant
measure carried by a closed geodesic is a semiclassical measure
arising from a sequence of eigenfunctions. Jakobson and Zelditch
then use the fact that the closed convex hull of such measures is
the set of all invariant measures. Following these ideas, the
result of Jakobson and Zelditch was extended in \cite{AzagraMacia}
to manifolds of constant positive sectional curvature~: any
invariant measure can be obtained as the semiclassical measure
arising from a sequence of eigenfunctions.

The spectrum of the Laplacian consists of clusters of bounded
width, centered at the points $\left(  k+\beta\right)  ^{2},$
$k\in\mathbb{Z}$, where $\beta>0$ is a constant depending on the
geometry of $M$ -- note that the spectrum of the sphere
$\mathbb{S}^{d}$ is exactly of this form. This was proved in
\cite{DuistermaatGuillemin, WeinsteinZoll, CdVSpectreZoll}, see
also \cite{UribeZelditch, ZelditchDegenerate, ZelditchZoll} for
more precise results on the structure of the spectrum. Using this
fact, it is possible to show that Jakobson and Zelditch's result
also holds for compact rank one symmetric spaces, see
\cite{MaciaZoll}. However, the statement for eigenfunctions of the
Laplacian on general Zoll manifolds does not seem to be known.

The situation for the time-dependent equation \eqref{e:densities}
is clearer. Consider a general sequence of initial data $\left(
u_{n}^{0}\right)  $ normalized in $L^{2}\left(  M\right)  $ and
chose $\left(  h_{n}\right)  $ as in Section
\ref{s:ltsemiclassical}. Suppose moreover that $\mu_{0}$ is the
unique semiclassical measure of this sequence and that $w_{h_{n}}$
converges to $\mu(t)$ as given by Theorem \ref{t:convergence}.
Maci\`{a} has proved in \cite{MaciaAv} the following result
relating $\mu_{0}$ to $\mu(t)$.

\begin{theorem}
\label{t:zoll}Let $\left(  M,g\right)  $ be a Zoll manifold and
$\mu_{0}$ and $\mu(t)$ be as above. Suppose $\mu_{0}\left( \left\{
\xi=0\right\}  \right) =0$. Then, for every $a\in C_{c}\left(
T^{\ast}M\right)  $ and a.e.
$t\in\mathbb{R}$ the following holds:%
\[
\int_{T^{\ast}M}a\left(  x,\xi\right)  \mu(t)\left( dx,d\xi\right)
=\int_{T^{\ast}M}\left\langle a\right\rangle \left( x,\xi\right)
\mu _{0}\left(  dx,d\xi\right)  ,
\]
where $\left\langle a\right\rangle $ is the average of $a$ along
the geodesic
flow.\footnote{That is:%
\[
\left\langle a\right\rangle \left(  x,\xi\right)
:=\lim_{T\rightarrow\infty }\frac{1}{T}\int_{0}^{T}a\circ
g_{s}\left(  x,\xi\right)  ds.
\]
}
\end{theorem}

Note, in particular, that $\mu(t)$ does not depend on $t$; if in
addition $\mu_{0}$ is an invariant measure then $\mu(t)=\mu_{0}$
for almost every $t\in\mathbb{R}$. When $\left(  u_{n}^{0}\right)
$ is a wave-packet (\ref{e:WP}) whose semiclassical measure is
$\delta_{\left(  x_{0},\xi
_{0}\right)  }$, Theorem \ref{t:zoll} implies that%
\[
\mu(t)=\lim_{T\rightarrow\infty}\frac{1}{T}\int_{0}^{T}\delta_{g^{s}\left(
x_{0},\xi_{0}\right)  }ds,
\]
in other words, $\mu(t)$ is the orbit measure on the geodesic
issued from $\left(  x_{0},\xi_{0}\right)  $. From the fact that
the closed convex hull of such measures is the whole set of
invariant measures, the following consequence is obtained.

\begin{corollary}
\label{c:allmeasures}Suppose $\left(  M,g\right)  $ is a Zoll
manifold. Then every invariant measure in $\mathcal{P}\left(
T^{\ast}M,g^{t}\right)  $ can be obtained as the semiclassical
measure (in the sense of Theorem \ref{t:convergence}) of some
sequence of initial data in $L^{2}\left( M\right)  $.
\end{corollary}

As mentioned in the introduction, this shows that Strichartz
estimates fail in Zoll manifolds. Combining Corollary
\ref{c:allmeasures} with Lebeau's result \cite{LebeauSchrod92}
gives the following.

\begin{corollary}
\cite{MaciaDispersion} Let $\left(  M,g\right)  $ be a Zoll
manifold, $T>0$ and $U\subset M$ an open set. Then condition
(\ref{e:GCC}) holds for $U$ if and only if the observability
property for the Schr\"{o}dinger flow holds for $U$ and $T$.
\end{corollary}

Note that in this result, $T>0$ can be chosen arbitrarily small,
since it does not play a role in condition (\ref{e:GCC}).

\section{The flat torus\label{s:torus}}
The geodesic flow on the cotangent bundle of the flat torus $\mathbb{T}%
^{d}:=\mathbb{R}^{d}/2\pi\mathbb{Z}^{d}$ is the prototype of a
non-degenerate
completely integrable Hamiltonian system. It has a simple explicit expression:%
\[
g^{s}\left(  x,\xi\right)  =\left(  x+s\xi,\xi\right)  .
\]
For each $\xi\in\mathbb{R}^{d}$ the torus
$\mathbb{T}^{d}\times\left\{
\xi\right\}  $ is an invariant Lagrangian submanifold of $T^{\ast}%
\mathbb{T}^{d}\cong\mathbb{T}^{d}\times\mathbb{R}^{d}$. In each of
these tori, the dynamics of the geodesic flow can be described in
terms of the order of resonance of $\xi$. More precisely, consider
the following primitive submodule of $\mathbb{Z}^{d}$,\footnote{A
submodule $\Lambda\subset\mathbb{Z}^{d}$ is \emph{primitive} if it
equals the intersection of $\mathbb{Z}^{d}$ with its
linear span $\left\langle \Lambda\right\rangle $ over $\mathbb{R}$.}%
\[
\Lambda_{\xi}:=\left\{  k\in\mathbb{Z}^{d}\;:\;k\cdot\xi=0\right\}
.
\]
Then the orbit issued from $\left(  x,\xi\right)  $, for any $x\in
\mathbb{T}^{d}$, is dense in a torus of dimension $d-\operatorname*{rk}%
\Lambda_{\xi}$; this quantity is sometimes called the order of
resonance of $\xi$. In particular, such a trajectory is periodic
(and non-constant) if $\operatorname*{rk}\Lambda_{\xi}=d-1$, and
dense on $\mathbb{T}^{d}$ when $\Lambda_{\xi}=\left\{  0\right\}
$.

If $\operatorname*{rk}\Lambda_{\xi}>0$ then $\xi$ is said to be
\emph{resonant}; this means that there exists $k\in\mathbb{Z}^{d}%
\setminus\left\{  0\right\}  $ such that $k\cdot\xi=0$. We shall
denote by $\Omega$ the set of all $\xi\in\mathbb{R}^{d}$ that are
resonant; they play an important role in the results we present
below.

Let us first recall some existing results for the case in which
the sequence of initial data consists of eigenfunctions of the
Laplacian. Let $\left( u_{n}\right)  $ be such that $-\Delta
u_{n}=\lambda_{n}u_{n}$, with $\left\Vert u_{n}\right\Vert
_{L^{2}\left(  \mathbb{T}^{d}\right)  }=1$ and
$\lambda_{n}\rightarrow\infty$.\ Clearly, $\lambda_{n}=\left\vert
k_{n}\right\vert ^{2}$ for some $k_{n}\in\mathbb{Z}^{d}$ and the
corresponding eigenfunction $u_{n}$ is a linear combination of
exponentials $e^{ik\cdot x}$ with $\left\vert k\right\vert
=\left\vert k_{n}\right\vert $. When $d=1$, the multiplicity of
$\lambda_{n}>0$ is equal to two and it follows that the weak
limits of the densities $\left\vert u_{n}\right\vert ^{2}$ are
constant. As soon as $d\geq2$ the multiplicity of $\lambda_{n}$
tends to infinity as $n\rightarrow\infty$, and the structure of
the limits becomes less evident. The following inequality is due
to Cooke \cite{CookeTorus2d} and Zygmund \cite{Zygmund74}: there
exists $C>0$ such that if $u$ is an eigenfunction of
the Laplacian on $\mathbb{T}^{2}$ then%
\begin{equation}
\left\Vert u\right\Vert _{L^{4}\left(  \mathbb{T}^{2}\right) }\leq
C\left\Vert u\right\Vert _{L^{2}\left( \mathbb{T}^{2}\right)  }.
\label{e:zygmund}%
\end{equation}
In particular, this implies that any accumulation point $\nu$ of
$\left\vert u_{n}\right\vert ^{2}$, in the weak-$\ast$ topology of
$\mathcal{M}_{+}\left( \mathbb{T}^{2}\right)  $, is in
$L^{2}\left(  \mathbb{T}^{2}\right)  $. This result was greatly
improved by Jakobson \cite{JakobsonTori97}, who showed that $\nu$
is in fact a trigonometric polynomial whose frequencies lie in at
most two circles centered at the origin. It is not known whether
an estimate such as (\ref{e:zygmund}) holds when $d\geq3$ (for
frequency dependent estimates see \cite{BourgainEigenf93}).
However, Bourgain has proved that any limit measure $\nu$ is
absolutely continuous with respect to the Lebesgue measure; in
fact, $\nu$ has additional regularity: if $\nu\left(  x\right)
=\sum _{k\in\mathbb{Z}^{d}}c_{k}e^{ik\cdot x}$ then
$\sum_{k\in}\left\vert c_{k}\right\vert ^{d-2}<\infty$. These
results are proved in \cite{JakobsonTori97, JakNadToth01,
Aissiou}, and rely on a deep understanding of the geometry of
lattice points in $\mathbb{R}^{d}$. The proof does not use
semiclassical analysis nor the relation between the Schr\"odinger
equation and the geodesic flow.

When the sequence of initial data $\left(  u_{n}^{0}\right)  $ is
not formed by eigenfunctions much less is known. The analogue of
(\ref{e:zygmund}) is
this setting is the estimate for $p=4$ and $d=1$:%
\[
\left(  \int_{0}^{1}\left\Vert e^{it\Delta/2}u\right\Vert
_{L^{p}\left( \mathbb{T}^{d}\right)  }^{p}dt\right)  ^{1/p}\leq
C\left\Vert u\right\Vert _{L^{2}\left(  \mathbb{T}^{d}\right)  }.
\]
However, no such inequality is known to hold when $d\geq2$
(Bourgain has made some conjectures in that direction
\cite{BourgainStrichartz93}).

The following holds.

\begin{theorem}
\label{t:abscon}Let $\left(  u_{n}^{0}\right)  $ be a bounded
sequence in $L^{2}\left(  \mathbb{T}^{d}\right)  $. If $\nu\in
L^{\infty}\left( \mathbb{R};\linebreak\mathcal{M}_{+}\left(
\mathbb{T}^{d}\right)  \right)  $
is obtained as the weak limit of $\left\vert e^{it\Delta/2}u_{n}%
^{0}\right\vert ^{2}$ (in the sense of (\ref{e:wL})) then $\nu$ is
absolutely continuous with respect to Lebesgue measure.
\end{theorem}

This result was proved by Bourgain \cite{BourgainQL97} using fine
results on the distribution of lattice points on paraboloids. It
can be also deduced as a consequence of the results of Maci\`{a}
\cite{MaciaTorus} and Anantharaman and Maci\`{a}
\cite{AnantharamanMacia} which we describe below. In
\cite{AnantharamanMacia} it is shown that Theorem \ref{t:abscon}
also holds for more general Hamiltonians of the form
$\frac{1}{2}\Delta + V(t,x)$. In the case of a coherent state \eqref{e:propatorus} or of a Lagrangian state on a torus,  the explicit computations of the densities
\eqref{e:densities}  and of their limits \eqref{e:wL} are presented in Propositions 13 and 14 of \cite{MaciaTorus}.

The proof of Theorem \ref{t:abscon} given in \cite{AnantharamanMacia,
MaciaTorus} relies on the structure of the geodesic flow on the
torus and is better understood in terms of semiclassical measures.

Suppose $\mu_{0}$ is a semiclassical measure of $\left(
u_{n}^{0}\right)  $
and that $\mu\in$\linebreak$L^{\infty}\left(  \mathbb{R};\mathcal{M}%
_{+}\left(  T^{\ast}\mathbb{T}^{d}\right)  \right)  $ is a
time-dependent semiclassical measure for $\left(
e^{it\Delta/2}u_{n}^{0}\right)  $, obtained as a weak-$\ast$ limit
as in Theorem \ref{t:convergence}. When $\mu_{0}\left(
\mathbb{T}^{d}\times\Omega\right)  =0$ (where $\Omega$ is the set
of resonant vectors defined above), it has been shown in
\cite{MaciaAv} that for almost
every $t\in\mathbb{R}$:%
\[
\mu(t)=\frac{1}{\left(  2\pi\right)  ^{d}}dx\otimes\int_{\mathbb{T}^{d}}%
\mu_{0}\left(  dy,\cdot\right)  .
\]
This can be seen as an analogue of Theorem \ref{t:zoll} in this
context, since for $\xi\in\mathbb{R}^{d}\setminus\Omega$ and $a\in
C_{c}\left(  T^{\ast
}\mathbb{T}^{d}\right)  $ one has for every $x\in\mathbb{T}^{d}$:%
\[
\lim_{T\rightarrow\infty}\frac{1}{T}\int_{0}^{T}a\circ g^{s}\left(
x,\xi\right)  ds=\frac{1}{\left(  2\pi\right)  ^{d}}\int_{\mathbb{T}^{d}%
}a\left(  y,\xi\right)  dy.
\]
Therefore, the non-trivial part of the structure of $\mu(t)$ is
that corresponding to its restriction to
$\mathbb{T}^{d}\times\Omega$. A first insight on the complexity of
this restriction is provided by a construction in
\cite{MaciaAv}. Two sequences of initial data exist such that both are $h_{n}%
$-oscillating for a common scale $\left(  h_{n}\right)  $ and have
as
semiclassical measure:%
\[
\mu_{0}\left(  x,\xi\right)  =\left\vert \rho\left(  x\right)
\right\vert ^{2}dx\delta_{\xi_{0}}\left(  \xi\right)
\]
with $\xi_{0}\in\Omega$ and $\rho\in L^{2}\left(
\mathbb{T}^{d}\right)  $ with $\left\Vert \rho\right\Vert
_{L^{2}\left(  \mathbb{T}^{d}\right)  }=1$. However, the
time-dependent measures corresponding to their orbits by the
Schr\"{o}dinger flow differ, and are respectively equal to:%
\begin{equation}
\left\vert e^{it\Delta/2}\rho\right\vert ^{2}dx\otimes\delta_{\xi_{0}}%
,\quad\text{and}\quad\frac{1}{\left(  2\pi\right)
^{d}}dx\otimes\delta
_{\xi_{0}}. \label{e:differentsm}%
\end{equation}
Therefore, in contrast with the situation in Theorem \ref{t:zoll},
the semiclassical measure of the initial data no longer determines
that of the evolutions. Moreover, those limiting measures $\mu(t)$
may have a non-trivial dependence on $t$.

A precise formula relating $\mu(t)$ to the sequence of initial
data is presented in \cite{MaciaTorus} for the two-dimensional
case and, more generally, in \cite{AnantharamanMacia} for any
dimension. For the sake of simplicity, let us give here the result
corresponding to $d=2$.

Start noticing that the set of resonant directions can be written
as a
disjoint union%
\begin{equation}
\Omega=%
%TCIMACRO{\dbigsqcup _{\omega\in\mathcal{R}}}%
%BeginExpansion
{\displaystyle\bigsqcup_{\omega\in\mathcal{R}}}
%EndExpansion
R_{\omega}\sqcup\left\{  0\right\}  , \label{e:resonant}%
\end{equation}
where $\mathcal{R}$ is formed by those vectors in $\mathbb{Z}^{2}$
whose components are relatively prime integers and
$R_{\omega}:=\omega^{\perp }\setminus\left\{  0\right\}  $. If a
measure $\mu\in\mathcal{P}\left(
T^{\ast}\mathbb{T}^{2},g^{s}\right)  $ is invariant, its
restriction to $\mathbb{T}^{2}\times\left(
\mathbb{R}^{2}\setminus\Omega\right)  $ and $\mathbb{T}^{2}\times
R_{\omega}$ enjoy more regularity on the $x$-variable than it is
expected \emph{a priori}. In fact, the restriction of $\mu$ to
these sets is constant with respect to $x$ along certain
directions \cite{AnantharamanMacia}: given $v\in\mathbb{R}^{2}$
write $\tau_{v}\left( x,\xi\right)  :=\left(  x+v,\xi\right)  $;
then $\left(  \tau_{v}\right)
_{\ast}\mu\rceil_{\mathbb{T}^{2}\times\left(
\mathbb{R}^{2}\setminus
\Omega\right)  }=\mu\rceil_{\mathbb{T}^{2}\times\left(  \mathbb{R}%
^{2}\setminus\Omega\right)  }$ for all $v\in\mathbb{R}^{2}$ (and
therefore, $\mu\rceil_{\mathbb{T}^{2}\times\left(
\mathbb{R}^{2}\setminus\Omega\right) }$ is constant in $x$),
whereas $\left(  \tau_{v}\right)  _{\ast}\mu
\rceil_{\mathbb{T}^{2}\times
R_{\omega}}=\mu\rceil_{\mathbb{T}^{2}\times R_{\omega}}$ holds for
every $v\in\omega^{\perp}$.

If $\mu(t)$ is a time-dependent semiclassical measure (as given by
Theorem \ref{t:convergence}) of a sequence $\left(
e^{it\Delta/2}u_{n}^{0}\right)  $ then it turns out that, besides
from the fact that $\mu(t)$ is invariant,
$\mu(t)\rceil_{\mathbb{T}^{2}\times R_{\omega}}$ enjoys additional
regularity in the directions in ${\mathbb{R}} \omega$. The reason
for this is that time-averaging produces a second
microlocalization around the lines $\omega^{\perp}$ which neglects
the contribution of the fraction of the energy of $\left(
e^{it\Delta/2}u_{n}^{0}\right)  $ that goes to infinity in the
direction $\omega$. In other words, if $a\in C^{\infty}\left(
\mathbb{T}^{2}\right)  $ is a function whose non vanishing Fourier
modes correspond to frequencies in $\omega$ and $\varphi\in
C_{c}^{\infty}\left(  \mathbb{R}^{2}\right)  $
vanishes in a neighborhood of $\omega^{\perp}$ then one has (see \cite{MaciaTorus}):%
\[
\int_{a}^{b}e^{-it\Delta/2}\operatorname*{Op}\nolimits_{h}\left(
a\otimes\varphi\right)  e^{it\Delta/2}dt=\left( \frac{b-a}{\left(
2\pi\right) ^{2}}\int_{\mathbb{T}^{2}}adx\right)  \varphi\left(
hD_{x}\right) +\mathcal{O}\left(  h\right)  .
\]

Denote by $L_{\omega}^{p}\left(  \mathbb{T}^{2}\right)  $ the
space of functions $a\in L^{p}\left(  \mathbb{T}^{2}\right)  $
such that $a\circ \tau_{v}=a$ for $v\in\omega^{\perp}$. The
following result holds (see \cite{MaciaTorus, AnantharamanMacia}).

\begin{theorem}
\label{t:2d}Let $\left(  u_{n}^{0}\right)  $ be a sequence
normalized in $L^{2}\left(  \mathbb{T}^{2}\right)  $. Suppose that
$\mu(t)$ is a semiclassical measure of $\left(
e^{it\Delta/2}u_{n}^{0}\right)  $, in the sense of Theorem
\ref{t:convergence}. Then for every $\omega\in\mathcal{R}$ there
exists a measure $\rho_{\omega}$, defined on $R_{\omega}$ and
taking values in the space of trace-class operators on
$L_{\omega}^{2}\left( \mathbb{T}^{2}\right)  $ such that for
a.e.\emph{ }$t\in\mathbb{R}$, every $a\in
L_{\omega}^{\infty}\left(  \mathbb{T}^{2}\right)  $ and every
$\varphi\in C_{c}\left(  \mathbb{R}^{2}\right)  $ we have:%
\begin{equation}
\int_{\mathbb{T}^{2}\times R_{\omega}}a\left(  x\right)
\varphi\left( \xi\right)  \mu(t)\left(  dx,d\xi\right)
=\int_{R_{\omega}}\varphi\left( \xi\right) \operatorname{tr}\left(
m_{a}e^{-it\Delta/2}\rho_{\omega}\left(
d\xi\right)  e^{it\Delta/2}\right)  , \label{e:traceclass}%
\end{equation}
where $m_{a}$ denotes the operator acting by multiplication by $a$
in $L_{\omega}^{2}\left(  \mathbb{T}^{2}\right)  $.
\end{theorem}

>From this, it follows that $\mu(t)\rceil_{\mathbb{T}^{2}\times
R_{\omega}}$ is absolutely continuous with respect to the
$x$-variable, and because of (\ref{e:resonant}), that
$\mu(t)\rceil_{\mathbb{T}^{2}\times\left(
\mathbb{R}^{2}\setminus\left\{  0\right\}  \right)  }$ is
absolutely
continuous. It is also possible to prove that $\mu(t)\rceil_{\mathbb{T}%
^{2}\times\left\{  0\right\}  }$ is also given by a formula
similar to (\ref{e:traceclass}), from which Theorem \ref{t:abscon}
follows. The ``measures'' $\rho_{\omega}$ depend only on the
(sub)sequence $\left( u_{n^{\prime}}^{0}\right)  $, however they
are not determined by the semiclassical measure of the initial
data alone, and are responsible for the phenomena presented in
example (\ref{e:differentsm}).

The generalization of Theorem \ref{t:2d} to dimensions higher than
two is non-trivial. The main reason for that is that in the
general case there is not a decomposition of the set of resonant
frequencies as simple as
(\ref{e:resonant}). Therefore, the analogues of identities (\ref{e:traceclass}%
) in this case are obtained by an iterative procedure that
requires to perform successive two-microlocalizations along nested
sequences of linear subspaces contained in the resonant set
$\Omega$, see \cite{AnantharamanMacia}.

As pointed out in the introduction, Theorem \ref{t:abscon} cannot
be used to obtain counterexamples to the validity of Strichartz
estimates for the Schr\"{o}dinger flow. On the other hand, it was
used in \cite{MaciaDispersion} to obtain an alternative proof of
Jaffard's result \cite{JaffardPlaques} on the observability of the
Schr\"{o}dinger flow on the bidimensional flat torus described in
Section \ref{s:observability}.

\section{Negatively curved manifolds\label{s:negcurv}}
In the case of negatively curved compact Riemannian manifolds, we
can make two contradictory remarks. The fact that the geodesic
flow has well-understood chaotic properties (to be precise, has
the {\emph{Anosov property}}) makes one very optimistic about the
good dispersive properties of the Schr\"odinger flow. This
motivates some very strong conjectures, such as the quantum unique
ergodicity conjecture (QUE) described below. On the other hand,
these same chaotic properties make it difficult to approximate the
Schr\"odinger dynamics by the geodesic dynamics~: the
quantum-classical correspondence is only valid for a relatively
short range of time (the Ehrenfest time), and this leaves little
hope to use it to prove those conjectures.

We first state the Snirelman theorem, whose proof can be found in
\cite{Sni, Zel87, CdV85}. On a smooth compact Riemannian manifold
$(M,g)$, take an orthonormal basis $(u_{n})$ of $L^{2}(M)$, formed
of eigenfunctions of the Laplacian ($-\Delta
u_{n}=\lambda_{n}u_{n}$, and $\lambda_{n}\nearrow\infty$). Assume
that the geodesic flow $g^{s}$ is ergodic with respect to the
Liouville measure. Write $h_{n}:=\lambda_{n}^{-1/2}$ and let
$w_{h_{n}}$ denote the
Wigner distribution of $u_{n}$. Then, there exists a subset ${\mathcal{S}%
}\subset{\mathbb{N}}$, of \emph{density 1}, such that the sequence $(w_{h_{n}%
})_{n\in{\mathcal{S}}}$ converges weakly to the Liouville measure.
Thus, the result says that a typical sequence of eigenfunctions
becomes equidistributed, both in the
\textquotedblleft$x$-variable\textquotedblright\ and in the
\textquotedblleft$\xi$-variable\textquotedblright. At this level
of generality, it is not well understood if the whole sequence
converges, or if there can be exceptional subsequences with a
different limiting behavior. There are manifolds (or Euclidean
domains) with ergodic geodesic flows, but with exceptional
subsequences of eigenfunctions \cite{Hass10}. But these examples
have only been found very recently, and we stress the fact that
the proof is not constructive; \emph{a fortiori}, the exceptional
subsequences, whose existence is proved, are not exhibited
explicitly. Thus, one cannot say that the phenomenon is fully
understood.

The statement of the Snirelman theorem can be adapted to solutions
of the time-dependent Schr\"{o}dinger equation \cite{AnRiv}~: take
a sequence of initial conditions $(u_{n}^{0})$ chosen randomly
from a \textquotedblleft generalized orthonormal
family\textquotedblright, with characteristic lengths of
oscillations $h_{n}$ going to $0$. Denote by $w_{h_{n}}(t)$ the
Wigner distribution associated with $(e^{it\Delta/2}u_{n}^{0})$.
Then the sequence $\int_{0}^{1}w_{h_{n}}(t)dt$ converges to the
Liouville measure, in the probabilistic sense (the paper
\cite{AnRiv} provides a detailed statement, and a rate of
convergence for negatively curved manifolds).

Negatively curved manifolds have ergodic geodesic flows, but
actually the understanding of the chaotic properties of the flow
is so good that one could hope to go beyond the Snirelman theorem.
It may seem surprising that the question is still widely open,
even in the case of manifolds of constant negative curvature
(where the local geometry is completely explicit). The QUE
conjecture was stated by Rudnick and Sarnak \cite{SarnakSchur,
RudSar} for eigenfunctions of the Laplacian on a negatively curved
compact manifold. If $(u_{n})$ is a sequence of eigenfunctions of
the Laplacian ($-h_{n}^{2}\Delta u_{n}=u_{n}$ with
$h_{n}\longrightarrow0$) and $(w_{h_{n}})$ the associated Wigner
distributions, the conjecture says that $(w_{h_{n}})$ converges to
the Liouville measure. In other words, there are no exceptional
subsequences of eigenfunctions for which $(w_{h_{n}})$ converges
to an invariant measure other than Liouville. So far, the only
complete result is due to E. Lindenstrauss \cite{BouLin,
LindenQUE}, who proved the conjecture in the case where $M$ is an
\emph{arithmetic congruence surface}, and the eigenfunctions
$(u_{n})$ are common eigenfunctions of $\Delta$ and of the
\emph{Hecke operators}. Unfortunately, his proof relies a lot on
the use of the Hecke operators, and cannot be adapted to more
general situations.

There is a partial result, due to N. Anantharaman and S.
Nonnenmacher, which holds in great generality, on any compact
negatively curved manifold \cite{NAQUE, AnNon, Riv}. Let again
$(u_{n})$ be a sequence of eigenfunctions, and $\mu$ be a limit
point of the corresponding Wigner distributions $(w_{h_{n}})$. The
following result deals with the \emph{Kolmogorov-Sinai entropy }of
the invariant measure $\mu$.

The Kolmogorov-Sinai entropy is a functional
$h_{KS}~:\mathcal{P}\left( S^{\ast}M,g^{s}\right)
\longrightarrow{\mathbb{R}}_{+}$, from the set of
$g^{s}$-invariant probability measures to ${\mathbb{R}}_{+}$. The
shortest (though not always the most convenient) definition of the
entropy results from a theorem due to Brin and Katok \cite{BK83}.
For any time $T>0$, introduce a distance on $S^{\ast}M$,
\[
d_{T}(\rho,\rho^{\prime})=\max_{t\in\lbrack-T/2,T/2]}d(g^{t}\rho,g^{t}%
\rho^{\prime}),
\]
where $d$ is the distance built from the Riemannian metric. For
$\epsilon>0$, denote by $B_{T}(\rho,\epsilon)$ the ball of center
$\rho$ and radius $\epsilon$ for the distance $d_{T}$. When
$\epsilon$ is fixed and $T$ goes to infinity, it looks like a
thinner and thinner tubular neighborhood of the geodesic segment
$[g^{-\epsilon}\rho,g^{+\epsilon}\rho]$.

Let $\mu$ be a $g^{s}$--invariant probability measure on
$T^{\ast}M$. Then, for $\mu$-almost every $\rho\in T^{\ast}M$, the
limit
\begin{multline*}
\lim_{\epsilon\longrightarrow0}\,\liminf_{T\longrightarrow+\infty}-\frac{1}%
{T}\log\mu\big(B_{T}(\rho,\epsilon)\big)\\
=\lim_{\epsilon\longrightarrow0}\,\limsup_{T\longrightarrow+\infty}-\frac
{1}{T}\log\mu\big(B_{T}(\rho,\epsilon)\big)=:h_{KS}(\mu,\rho)
\end{multline*}
exists and it is called the local entropy of the measure $\mu$ at
the point $\rho$ (it is independent of $\rho$ if $\mu$ is
ergodic). The Kolmogo\-rov-Sinai entropy is the average of the
local entropies:
\[
h_{KS}(\mu)=\int h_{KS}(\mu,\rho)d\mu(\rho).
\]

In the case when $\mu$ is obtained from a limit of Laplace
eigenfunctions, the result of Anantharaman-Nonnenmacher says that
$h_{KS}(\mu)>0$. This is a strong restriction, for instance, $\mu$
cannot be entirely concentrated on a countable union of closed
geodesics. The result has been extended to the time-dependent
context of Theorem \ref{t:convergence} by
Anantharaman--Rivi\`{e}re, who showed that, for a.e
$t\in\mathbb{R}$, $\mu(t)$ has positive entropy. In the case of
manifolds of constant curvature $-1$, and dimension $d$, there is
actually an explicit lower bound, for the eigenfunction case
\cite{AnNon} and more generally for the time-dependent case
\cite{AnRiv}. In the latter case, one can disintegrate
$\mu(t)(x,\xi)=\int \mu_{E}(t)(x,\xi)\nu(dE)$ where $\nu$ is a
positive measure, and $\mu_{E}(t)$ is a $g^{s}$-invariant
probability measure supported on the energy layer
$\{\Vert\xi\Vert_{x}^{2}=E\}$. Then, one has, $dt\otimes\nu$
almost everywhere,
\[
h_{KS}(\mu_{E}(t))\geq\frac{d-1}{2}\sqrt{E}.
\]
We note that $\sqrt{E}$ is the speed of the geodesics on
$\{\Vert\xi\Vert _{x}^{2}=E\}$, and that $(d-1)\sqrt{E}$ is the
maximal entropy for invariant measures carried by this set (the
maximum is achieved only for the Liouville measure). The statement
for eigenfunctions is similar, with $\mu$ concentrated on one
single energy layer.

One can use this result to improve the Geometric Control condition
(\ref{e:GCC}) of \S \ref{s:observability}~:

\begin{theorem}
\cite{AnRiv} \label{t:observability} Let $M$ be a compact
Riemannian manifold of dimension $d$ and constant curvature
$\equiv-1$. Let $a$ be a smooth
function on $M$, and define the closed $g^{s}$-invariant subset of $S^{\ast}%
M$,
\[
K_{a}=\{\rho\in S^{\ast}M,a(g^{s}(\rho))=0\,\,\forall
s\in{\mathbb{R}}\}.
\]
Assume that the topological entropy of $K_{a}$ is
$<\frac{d-1}{2}$. Then, for all $T>0$, there exists $C_{T,a}>0$
such that, for all $u$~:
\begin{equation}
\Vert u\Vert_{L^{2}(M)}^{2}\leq C_{T,a}\int_{0}^{T}\Vert
ae^{\imath t\Delta
/2}u\Vert_{L^{2}(M)}^{2}dt. \label{e:observability}%
\end{equation}

\end{theorem}

For the purposes of this survey, we simply define the topological
entropy of $K_{a}$ as the supremum of Kolmogorov-Sinai entropies,
for all invariant measures supported on $K_{a}$. For manifolds of
constant curvature $\equiv-1$, the topological entropy is closely
related to the more familiar notion of Hausdorff dimension, at
least if $K_{a}$ is locally maximal~: saying that the topological
entropy of $K_{a}$ is $<\frac{d-1}2$ is equivalent to $K_{a}$
having Hausdorff dimension $<d$ \cite{PesSad}.

We see that the Geometric Control condition ($K_{a}$ empty) is
weakened by only assuming that $K_{a}$ have small Hausdorff
dimension.
Examples of such functions $a$ on a negatively curved surface of genus $g$ are given in \cite{AnRiv}~: one takes a decomposition of the surface into $2g-2$ ``hyperbolic pairs of pants'' with very long boundary components, and takes $a$ to be non-zero in a neighbourhood of these $3g-3$ curves. It would be interesting if one could enrich this list of examples.

\bibliographystyle{amsalpha}
\def\cprime{$'$}
\providecommand{\bysame}{\leavevmode\hbox
to3em{\hrulefill}\thinspace}
\providecommand{\MR}{\relax\ifhmode\unskip\space\fi MR }
\providecommand{\MRhref}[2]{%
  \href{http://www.ams.org/mathscinet-getitem?mr=#1}{#2}
} \providecommand{\href}[2]{#2}

\end{document}